\documentclass{amsart}
\usepackage[utf8]{inputenc}
\usepackage{amsthm,amsmath,amssymb,tikz-cd,mathrsfs,enumitem}

\makeatletter
\@namedef{subjclassname@2010}{%
  \textup{2010} Mathematics Subject Classification}
\makeatother

\title[$\Sigma^*$-Morita equivalence]{Hilbert $C^*$-modules over $\Sigma^*$-algebras II: \\ $\Sigma^*$-Morita equivalence}
\author{Clifford A. Bearden}
\address{Department of Mathematics\\The University of Texas at Tyler\\Tyler, TX 75799, USA}
\email{cbearden@uttyler.edu}

\date{}

\newtheorem{thm}{Theorem}[section]
\newtheorem{prop}[thm]{Proposition}
\newtheorem{lemma}[thm]{Lemma}
\newtheorem{cor}[thm]{Corollary}

\theoremstyle{definition}
\newtheorem{defn}[thm]{Definition}

\newtheorem{note}[thm]{Note}
\newtheorem{piece}[thm]{}

\newtheorem{example}[thm]{Example}

\newtheorem{notation}[thm]{Notation}

\newcommand{\N}{\mathbb{N}}

\newcommand{\K}{\mathbb{K}}
\newcommand{\B}{\mathbb{B}}

\newcommand{\Hil}{\mathcal{H}}
\newcommand{\Kil}{\mathcal{K}}
\newcommand{\bw}{\Sigma^*}
\newcommand{\swb}{\Sigma^*_{\mathfrak{B}}\text{-}}

\newcommand{\ssto}{\xrightarrow{\sigma}}
\newcommand{\rmto}[1]{\xrightarrow{\sigma_{#1}}}
\newcommand{\lmto}[1]{\xrightarrow{{}_{#1}\sigma}}

\newcommand{\fA}{\mathfrak{A}}
\newcommand{\fB}{\mathfrak{B}}
\newcommand{\fC}{\mathfrak{C}}

\newcommand{\fX}{\mathfrak{X}}
\newcommand{\fY}{\mathfrak{Y}}

\newcommand{\Bw}{\mathscr{B}}

\newcommand{\bwtp}{\mathfrak{X} \otimes_{\mathfrak{B}}^{\Sigma^*} \mathfrak{Y}}

\newcommand{\blink}[1]{\mathcal{L}^\mathscr{B}(#1)}

\begin{document}

\renewcommand{\thefootnote}{}
%
%



\begin{abstract}
In previous work, we defined and studied \emph{$\bw$-modules}, a class of Hilbert $C^*$-modules over $\bw$-algebras (the latter are $C^*$-algebras that are sequentially closed in the weak operator topology). The present work continues this study by developing the appropriate $\bw$-algebraic analogue of the notion of strong Morita equivalence for $C^*$-algebras. We define strong $\Sigma^*$-Morita equivalence, prove a few characterizations, look at the relationship with equivalence of categories of a certain type of Hilbert space representation, study $\Sigma^*$-versions of the interior and exterior tensor products, and prove a $\Sigma^*$-version of the Brown-Green-Rieffel stable isomorphism theorem.
\end{abstract}

\maketitle

\section{Introduction}


Strong Morita equivalence for $C^*$-algebras is an equivalence relation coarser than $*$-isomorphism, but fine enough to preserve many distinguishing structures and properties a $C^*$-algebra can have (e.g.\ ideal structure, representation theory, and $K$-theory in the $\sigma$-unital case). Originally a concept in pure ring theory, Morita equivalence was imported into operator algebra theory by M.\ Rieffel in \cite{Rie,Rie2}, where he defined and initiated the study of $C^*$-algebraic Morita equivalence and the corresponding version for $W^*$-algebras. It has since taken a central role in operator algebra theory, and appears in a number of the most important classification results for $C^*$-algebras (e.g.\ the Kirchberg-Phillips classification for Kirchberg algebras \cite[Theorem~8.4.1]{Ror}, the Dixmier-Douady classification for continuous-trace $C^*$-algebras \cite[Theorem~5.29]{RW}, and the classification of unital graph $C^*$-algebras recently announced by Eilers, Restorff, Ruiz, and S\o rensen \cite{ERRS}). The main goal of the present work is to define and study a version of Morita equivalence for a special class of $C^*$-algebras called ``$\Sigma^*$-algebras".

First defined and studied by E. B. Davies in \cite{Dav68}, a (concrete) \emph{$\bw$-algebra} is a $C^*$-subalgebra of $B(\Hil)$ that is closed under limits of sequences converging in the weak operator topology (WOT). Evidently, every von Neumann algebra is a $\bw$-algebra, but the converse does not hold. Conceptually, the theory of $\bw$-algebras may be considered as a halfway world between general $C^*$-algebra theory and von Neumann algebra theory---indeed, many (but certainly not all) von Neumann algebraic concepts and techniques have some sort of analogue for $\bw$-algebras. In previous work \cite{Bear}, we defined and studied the $\bw$-algebraic analogue of $W^*$-modules (the latter are selfdual $C^*$-modules over $W^*$-algebras), which we called ``$\bw$-modules". These will be defined in the next section and used extensively throughout this paper. We refer to \cite{Bear} for more discussion about the motivation and underlying philosophy of this project.

The organization of this paper is as follows. In Section~2, we provide some background material on Morita equivalence for $C^*$-algebras and $W^*$-algebras (mostly as motivation for the analogues we develop later, but some of this material is used in this paper), and we cover relevant definitions, notations, and results for $\bw$-algebras and $\bw$-modules. In Section~3, we define strong $\bw$-Morita equivalence and $\bw$-imprivitivity bimodules, and prove a few basic results that are used throughout the rest of the paper. In Section~4, we give an account of the $\bw$-module interior tensor product of $\bw$-imprivitivity bimodules, with the main goal being to prove that $\bw$-Morita equivalence is transitive (hence an equivalence relation). Section~5 contains a discussion of ``weak Morita equivalence" for $C^*$-algebras and $\bw$-algebras (in the $W^*$-case, the analogous definition coincides with the ``strong" version). In Section~6, we prove a version of the ``full corners" characterization of strong Morita equivalence and some consequences. In Section~7, we develop the $\bw$-module exterior tensor product, which allows us in Section~8 to prove an analogue of the Brown-Green-Rieffel stable isomorphism theorem.

Some of the results in this paper were announced at the Workshop in Analysis and Probability at Texas A\&M University in July 2015 and the mini-workshop ``Operator Spaces and Noncommutative Geometry in Interaction" at Oberwolfach in February 2016. The author would like to thank the organizers of these conferences, as well as the NSF for funding during part of the time research was carried out for this paper.

As with \cite{Bear}, this work constitutes a significant part of the author's doctoral thesis at the University of Houston. The author would again like to express his deep gratitude for the guidance of his Ph.D. advisor, David Blecher.

\section{Background}

\subsection{$C^*$-modules, $W^*$-modules, and Morita equivalence}

Since the basic theory of Morita equivalence for $C^*$-algebras is well-known and covered in many texts, we will be brief here. We generally refer to \cite[Chapter~8]{BLM} for notation and results; other references include \cite{Lan, RW}, and \cite[Section~II.7]{Bl}.


Loosely speaking, a (right or left) module $X$ over a $C^*$-algebra $A$ is called a \emph{(right or left) $C^*$-module over $A$} if $X$ is equipped with an ``$A$-valued inner product" $\langle \cdot | \cdot \rangle : X \times X \to A$ (see any of the references above for the complete list of axioms, including, for example, $\langle x|x \rangle \geq 0$ for all $x \in X$) such that $X$ is complete in the canonical norm induced by this inner product. If $X$ is a right $C^*$-module, the inner product is taken to be linear and $A$-linear in the second variable and conjugate-linear in the first variable, and vice versa for left modules.

For clarity, we sometimes write the inner product with a subscript denoting its range (e.g.\ $\langle \cdot | \cdot \rangle_A$ for a right module and ${}_A \langle \cdot|\cdot \rangle$ for a left module). If $X$ is a right $C^*$-module, we will denote by $\overline X$ the \emph{adjoint module} of $X$ (often called the \emph{conjugate module}); see \cite{BLM}, last paragraph of 8.1.1.

If $X$ and $Y$ are two $C^*$-modules over $A$, $B_A(X,Y)$ will denote the Banach space of bounded $A$-module maps from $X$ to $Y$ with operator norm; $\mathbb{B}_A(X,Y)$ will denote the closed subspace of adjointable operators; and $\mathbb{K}_A(X,Y)$ will denote the closed subspace generated by operators of the form $|y \rangle \langle x| := y \langle x | \cdot \rangle$ for $y \in Y,$ $x \in X.$ If $X=Y$, the latter two of these spaces are $C^*$-algebras, and in this case, $X$ is a left $C^*$-module over both of these with inner product $|\cdot \rangle \langle \cdot|.$

Suppose $A \subseteq B(\Hil)$ is a $C^*$-algebra (assumed to be nondegenerately acting, although this is not strictly necessary for everything that follows) and that $X$ is a right $C^*$-module over $A$. We may consider $\Hil$ as a left module over $A$ and take the algebraic module tensor product $X \odot_A \Hil.$ This vector space admits an inner product determined by the formula $\langle x \otimes \zeta, y \otimes \eta \rangle = \langle \zeta, \langle x|y \rangle \eta \rangle$ for simple tensors (see \cite[Proposition~4.5]{Lan} for details), and we may complete $X \odot_A \Hil$ in the induced norm to yield a Hilbert space $X \otimes_A \Hil.$ (This is a special case of the interior tensor product of $C^*$-modules, cf.\ \ref{int tp constr} below.) Considering $A$ as a $C^*$-module over itself and taking the $C^*$-module direct sum $X \oplus A$, there is a canonical corner-preserving embedding of $B_A(X \oplus A)$ into $B((X \otimes_A \Hil) \oplus^2 \Hil)$ which allows us to concretely identify $X$, $A$, and many of the associated spaces of operators between $X$ and $A$ with spaces of Hilbert space operators between $\Hil$ and $X \otimes_A \Hil$ (see \cite[Proposition~2.1]{Bear} for more details). The following proposition records the part of this to which we will need to refer later:

\begin{prop} \label{mod embs}
If $X$ is a right (resp.\ left) $C^*$-module over a $C^*$-algebra $A \subseteq B(\Hil)$, then there is a canonical isometry $X \hookrightarrow B(\Hil, X \otimes_A \Hil)$  (resp.\ $X \hookrightarrow B(\overline X \otimes_A \Hil, \Hil)$).
\end{prop}

If $X$ is a left or right $C^*$-module over $A$, denote by $\langle X|X \rangle$ the linear span of the range of the $A$-valued inner product, i.e.
\[ \langle X|X \rangle := \mathrm{span}\{ \langle x|y \rangle : x,y \in X\}. \] Note that $\langle X|X \rangle$ is also the non-closed ideal in $A$ generated by $\{ \langle x|y \rangle : x,y \in X\}$. We say that $X$ is a \emph{full} $C^*$-module over $A$ if $\langle X|X \rangle$ is norm-dense in $A$.

A right $C^*$-module $X$ over a $A$ is called \emph{selfdual} if every bounded $A$-module map $X \to A$ is of the form $\langle x|\cdot \rangle$ for some $x \in X.$ A \emph{$W^*$-module} is a selfdual $C^*$-module over a $W^*$-algebra. (In our usage, a \emph{$W^*$-algebra} is an abstract $C^*$-algebra that admits a faithful representation as a von Neumann algebra.) A $W^*$-module $Y$ over a $W^*$-algebra $M$ is said to be \emph{$w^*$-full} if $\langle Y|Y \rangle$ is weak*-dense in $M$.

\begin{defn}
Two $C^*$-algebras $A$ and $B$ are said to be \emph{strongly $C^*$-Morita equivalent} if there exists an $A-B$ bimodule $X$ that is a full left $C^*$-module over $A$ and a full right $C^*$-module over $B$, and such that the two inner products and module actions on $X$ are related according to the following formula:
\[{}_A \langle x|y \rangle z = x \langle y|z \rangle_B \text{ for all } x,y,z \in X. \]
In this case $X$ is called an \emph{$A-B$ $C^*$-imprivitivity bimodule}.

Similarly, two $W^*$-algebras $M$ and $N$ are \emph{$W^*$-Morita equivalent} if there exists an $M-N$ bimodule $Y$ that is a $w^*$-full left $W^*$-module over $M$ and a $w^*$-full right $W^*$-module over $N$, and such that ${}_M \langle x|y \rangle z = x \langle y|z \rangle_N$ for all $x,y,z \in Y$. In this case $Y$ is called an \emph{$M-N$ $W^*$-imprivitivity bimodule}.
\end{defn}

\begin{piece}[Interior tensor product of $C^*$-modules] \label{int tp constr}
We sketch the construction of the interior tensor product of $C^*$-modules (see \cite[pgs.~38--41]{Lan} for details). Suppose $A$ and $B$ are $C^*$-algebras, $X$ is a right $C^*$-module over $A$, and $Y$ is a right $C^*$-module over $B$ such that there is a $*$-homomorphism $\rho: A \to \B_B(Y)$. The algebraic module tensor product $X \odot_A Y$ is then a right $B$-module which admits a $B$-valued inner product determined by the formula
\[ \langle x_1 \otimes y_1 | x_2 \otimes y_2 \rangle_B = \langle y_1 | \rho(\langle x_1 | x_2 \rangle_A) y_2 \rangle \]
for $x_1,x_2 \in X$ and $y_1,y_2 \in Y$. The completion of $X \odot_A Y$ in the norm induced by this inner product yields a right $C^*$-module $X \otimes_A Y$ over $B$, called the \emph{interior tensor product} of $X$ and $Y$.
\end{piece}

Now we record a couple of the most important results in basic $C^*$- and $W^*$-algebraic Morita equivalence theory. We will prove $\bw$-analogies for each of these in this work. Following the theorems is a note on some aspects of these theorems relevant to the present work, as well as explanations for some of the undefined terms.

\begin{thm} \label{Me er}
Strong $C^*$-Morita equivalence and $W^*$-Morita equivalence are equivalence relations strictly coarser than $*$-isomorphism.
\end{thm}

\begin{thm} \label{Me full corns}
Two $C^*$-algebras (resp.\ $W^*$-algebras) $A$ and $B$ are strongly $C^*$-Morita equivalent (resp.\ $W^*$-Morita equivalent) if and only if they are $*$-isomorphic to complementary full (resp.\ $w^*$-full) corners of a $C^*$-algebra (resp.\ $W^*$-algebra) $C$.
\end{thm}

\begin{thm} \label{BGRSIT}

\begin{enumerate}
\item[$\mathrm{(1)}$] Two $\sigma$-unital $C^*$-algebras $A$ and $B$ are strongly $C^*$-Morita equivalent if and only if $A$ and $B$ are stably isomorphic (i.e., $A \otimes \mathbb K \cong B \otimes \mathbb K$, where $\mathbb K$ is the space of compact operators on a separable Hilbert space).

\item[$\mathrm{(2)}$] Two $W^*$-algebras $M$ and $N$ are $W^*$-Morita equivalent if and only if there exists a Hilbert space $\Hil$ such that $M \overline{\otimes} B(\Hil) \cong N \overline{\otimes} B(\Hil)$ (where $\overline{\otimes}$ denotes the spatial von Neumann algebra tensor product).
\end{enumerate}
\end{thm}

\begin{note}
The $C^*$-algebra versions of all these are in the standard texts \cite{Lan, RW}. For the $W^*$-versions, see \cite[Section~8.5]{BLM}.

The only non-trivial part of Theorem~\ref{Me er} is transitivity, which is typically proved in the $C^*$-algebra case using the interior tensor product of $C^*$-modules outlined above. The proof of transitivity in the $W^*$-case may be done in many different ways, including via an analogue of the interior tensor product (see \cite[end of chapter notes for 8.5]{BLM} for a different way). In this paper, we develop a $\bw$-analogue of the interior tensor product to prove the corresponding result for strong $\bw$-Morita equivalence.

To explain the undefined terms in Theorem~\ref{Me full corns}, \emph{complementary corners} of a $C^*$-algebra $C$ are $C^*$-subalgebras of the form $pCp$ and $(1-p)C(1-p)$ for a projection $p$ in the multiplier algebra $M(C)$. A corner of a $C^*$-algebra (resp.\ $W^*$-algebra) is \emph{full} (resp.\ $w^*$-full) if the closed (resp.\ $w^*$-closed) ideal it generates is all of $C$. One direction of Theorem~\ref{Me full corns} is proved using an important construction called the \emph{linking algebra}, and we prove our analogue with the obvious $\bw$-version of the linking algebra.


Part (1) of Theorem~\ref{BGRSIT} is the incredible Brown-Green-Rieffel stable isomorphism theorem. A $C^*$-algebra is $\sigma$-unital if and only if it has a sequential contractive approximate identity, so in particular every separable $C^*$-algebra is $\sigma$-unital. There are examples of pairs of $C^*$-algebras that are strongly $C^*$-Morita equivalent but not stably isomorphic---the easiest example is $\mathbb C$ and $\mathbb K(\Hil)$ for non-separable $\Hil$ (see \cite[II.7.6.10]{Bl} for mention of a fancier and probably more satisfying example). Note the interesting fact that in Theorem~\ref{BGRSIT} (2), there are no restrictions on the $W^*$-algebras involved.
\end{note}

\subsection{$\bw$-algebras and $\bw$-modules}

We review now some of the relevant definitions and notations for $\bw$-algebras and $\bw$-modules. Compared to \cite{Bear}, we take a more abstract perspective in this paper since we will sometimes have to deal simultaneously with more than one representation of a given $\bw$-algebra.

\begin{defn}[\cite{Dav68}]
A \emph{concrete $\bw$-algebra} is a nondegenerate $C^*$-algebra $\fB \subseteq B(\Hil)$ that is closed under limits of WOT-convergent sequences; i.e., whenever $(b_n)$ is a sequence in $\fB$ that converges in the weak operator topology of $B(\Hil)$ to an operator $T$, then $T \in \fB.$

An \emph{abstract $\bw$-algebra} $(\fB, \mathscr S_\fB)$ is a $C^*$-algebra $\fB$ with a collection of pairs $\mathscr S_\fB = \{((b_n),b)\}$, each consisting of a sequence and an element in $\fB$, such that there exists a faithful representation $\pi: \fB \to B(\Hil)$ in which $\pi(\fB)$ is a concrete $\bw$-algebra and $((b_n),b) \in \mathscr S_\fB$ if and only if $\pi(b_n) \xrightarrow{WOT} \pi(b)$ in $B(\Hil)$. The sequences in $\mathscr S_\fB$ are called the \emph{$\sigma$-convergent} sequences of $\fB$, and we write $b_n \ssto b$ to mean $((b_n),b) \in \mathscr S_\fB$.

A \emph{$\bw$-representation} of an abstract $\bw$-algebra $(\fB, \mathscr S_\fB)$ is a nondegenerate representation $\pi: \fB \to B(\Hil)$ such that $b_n \ssto b$ implies $\pi(b_n) \xrightarrow{WOT} \pi(b)$. A \emph{faithful $\bw$-representation} is an isometric $\bw$-representation $\pi$ such that $\pi(\fB)$ is a concrete $\bw$-algebra and $b_n \ssto b$ if and only if $\pi(b_n) \xrightarrow{WOT} \pi(b)$.

Let $\fA$ and $\fB$ be abstract $\bw$-algebras. A map $\varphi: \fA \to \fB$ is \emph{$\sigma$-continuous} if $\varphi(a_n) \ssto \varphi(a)$ whenever $a_n \ssto a$ in $\fA$. A \emph{$\bw$-isomorphism} is a $*$-isomorphism $\psi: \fA \to \fB$ such that $\psi$ and $\psi^{-1}$ are $\sigma$-continuous. A \emph{$\bw$-subalgebra} of $\fB$ is a $C^*$-subalgebra $\fC$ closed under limits of $\sigma$-convergent sequences in $\fC$. A \emph{$\bw$-embedding} of $\fA$ into $\fB$ is an isometric $*$-homomorphism $\rho: \fA \hookrightarrow \fB$ such that $\rho(\fA)$ is a $\bw$-subalgebra of $\fB$ and $\rho$ is a $\bw$-isomorphism onto $\rho(\fA)$.
\end{defn}

\begin{note}
The phrase ``$\fB$ is a $\bw$-algebra" should be taken to mean that $\fB$ is an abstract $\bw$-algebra with an implicit collection of $\sigma$-convergent sequences; the phrase ``$\fB \subseteq B(\Hil)$ is a $\bw$-algebra" should be taken to mean that $\fB$ is a concrete $\bw$-algebra in $B(\Hil)$.

E. B. Davies in \cite{Dav68} provided a characterization for when a $C^*$-algebra $\fB$ with a collection of pairs $\{((b_n),b)\} \subseteq \ell^\infty(\fB) \times \fB$ is an abstract $\bw$-algebra.
\end{note}


\begin{defn}[\cite{Bear}] \label{ss-module def}
A right (resp.\ left) $C^*$-module $\fX$ over a concrete $\bw$-algebra $\fB \subseteq B(\Hil)$ is called a \emph{right (resp.\ left) $\bw$-module over $\fB$} if $\fX$ is WOT sequentially closed in $B(\Hil, \fX \otimes_\fB \Hil)$ (resp.\ in $B(\overline{\fX} \otimes_\fB \Hil, \Hil)$).

A right (resp.\ left) $C^*$-module $\fX$ over an abstract $\bw$-algebra $(\fB,\mathscr S_\fB)$ is called a \emph{right (resp.\ left) $\bw$-module over $\fB$} if the following holds: whenever $(x_n)$ is a sequence in $\fX$ such that $(\langle x_n| y \rangle)$ is $\sigma$-convergent for all $y \in \fX$, then there is an $x \in \fX$ such that $\langle x_n | y \rangle \ssto \langle x|y \rangle$.
\end{defn}

\begin{note}
The first definition above should be interpreted in reference to the embeddings in Proposition~\ref{mod embs}. It was the original definition given in \cite{Bear}.

The second definition is the abstract version, and is essentially equivalent. Indeed, by \cite[Proposition~3.4]{Bear}, a $\bw$-module over a concrete $\bw$-algebra $\fB$ is also $\bw$-module over $\fB$ considered as an abstract $\bw$-algebra; and conversely a $\bw$-module over an abstract $\bw$-algebra $\fB$ is a $\bw$-module over the concrete $\bw$-algebra $\pi(\fB)$ for any faithful $\bw$-representation $\pi$.
\end{note}

\begin{notation}
If $S$ is any subset of $B(\Hil)$, let $\mathscr B(S)$ (or $\mathscr B_\Hil(S)$ if there is a chance of ambiguity) denote the smallest WOT sequentially closed subset of $B(\Hil)$ containing $S$ (in general, $\Bw(S)$ is more than the set of limits of WOT-convergent sequences from $S$; e.g., $\Bw(S)$ contains limits of sequences consisting of limits of sequences from $S$). It is a simple exercise to show that $\Bw(A)$ is a $C^*$-algebra when $A$ is a $*$-subalgebra of $B(\Hil)$. Similarly, for a subset $S$ of a $\bw$-algebra $\fB$, denote by $\mathscr B(S)$ the smallest subset of $\fB$ that is closed under $\sigma$-convergent sequences and contains $S$. If $\Bw(S)=\fB$, we say that $S$ is \emph{$\sigma$-dense} in $\fB$, and if $\Bw(S)=S$, we say that $S$ is \emph{$\sigma$-closed in $\fB$}.

If $\fX$ is a right (resp.\ left) $\bw$-module over a $\bw$-algebra $\fB$ and $(x_n),x \in \fX$, we write
\[ x_n \rmto{\fB} x \quad (\text{resp. } x_n \lmto{\fB} x) \]
to mean that $\langle x_n|y \rangle \ssto \langle x|y \rangle$ for all $y \in \fX$, and we call this convergence \emph{$\sigma_\fB$-convergence} (resp.\ \emph{${}_\fB \sigma$-convergence}). We leave it as an exercise using Lemma~\ref{simplem} below to check that $x_n \rmto{\fB} x$ (resp.\ $x_n \lmto{\fB} x$) if and only if $x_n \xrightarrow{WOT} x$ in $B(\Hil, \fX \otimes_\fB \Hil)$ (resp.\ in $B(\overline{\fX} \otimes_\fB \Hil, \Hil)$) for any (hence every) faithful $\bw$-representation $\fB \hookrightarrow B(\Hil)$.
\end{notation}

Finally, we record some simple principles that we will use repeatedly:


\begin{lemma} \label{simplem}
Let $\Hil_1$ and $\Hil_2$ be Hilbert spaces.
\begin{enumerate}
\item[$\mathrm{(1)}$] A WOT-convergent sequence in $B(\Hil_1,\Hil_2)$ is bounded.

\item[$\mathrm{(2)}$] Let $\mathscr{T}_1, \mathscr{T}_2$ be sets of total vectors in $\Hil_1, \Hil_2$ respectively (i.e., the span of $\mathscr{T}_i$ is dense in $\Hil_i$).
A sequence $(T_n)$ in $B(\Hil_1,\Hil_2)$ is WOT-convergent if and only if $(T_n)$ is bounded and $\langle T_n \zeta, \eta \rangle$ converges for all $\zeta \in \mathscr{T}_1$ and $\eta \in \mathscr{T}_2$. Also, $T_n \xrightarrow{WOT} T$ if and only if $(T_n)$ is bounded and $\langle T_n \zeta, \eta \rangle \to \langle T \zeta, \eta \rangle$ for all $\zeta \in \mathscr{T}_1$ and $\eta \in \mathscr{T}_2$.
\end{enumerate}
\end{lemma}

\proof
We leave the proof as an exercise using the principle of uniform boundedness for (1), and using (1), the triangle inequality, and the correspondence between operators and bounded sesquilinear forms for (2).
\endproof

\section{Strong $\bw$-Morita equivalence}

In this section, $\fA$ and $\fB$ are always taken to be $\bw$-algebras.

\begin{defn}
A $\bw$-module $\fX$ over $\fB$ is \emph{$\sigma$-full over $\fB$} if $\Bw(\langle \fX|\fX \rangle) = \fB$.
\end{defn}

\begin{defn} \label{Me}
An $\fA-\fB$ bimodule $\fX$ is called an \emph{$\fA-\fB$ $\bw$-imprivitivity bimodule} if:
\begin{enumerate}
\item $\fX$ is a $\sigma$-full left $\bw$-module over $\fA$ and a $\sigma$-full right $\bw$-module over $\fB$;
\item ${}_\fA \langle x|y \rangle z = x \langle y|z \rangle_\fB$ for all $x,y,z \in \fX$;
\item for a sequence $(x_n) \in \fX$ and $x \in \fX$, $x_n \lmto{\fA} x$ if and only if $x_n \rmto{\fB} x$.
\end{enumerate}
If there exists an $\fA-\fB$ $\bw$-imprivitivity bimodule, we say that $\fA$ and $\fB$ are \emph{strongly $\bw$-Morita equivalent}.
\end{defn}

\begin{note} \label{Me note}
In the $C^*$- and $W^*$-settings, the analogue of condition (3) above is automatic, but we do not know if (3) is automatic in our case. One way to show that (3) does not always hold would be to exhibit two $\bw$-algebras that are isomorphic as $C^*$-algebras but are not $\bw$-isomorphic. (The existence of such algebras is an open question.) Indeed, suppose $\fA$ and $\fB$ are $\bw$-algebras that are $*$-isomorphic, but not $\bw$-isomorphic. If $\fX$ is the underlying $C^*$-algebra, then $\fX$ evidently satisfies conditions (1) and (2), but not (3).
\end{note}

\begin{example} \label{ssMe examps}
(Examples of $\bw$-imprivitivity bimodules.)
\begin{enumerate}
\item If $\fX$ is any right $\bw$-module over $\fB$, then by \cite[Lemma~3.7 and Theorem~3.8]{Bear}, $\fX$ is a $\Bw(\K_\fB(\fX))-\Bw(\langle \fX|\fX \rangle)$ $\bw$-imprivitivity bimodule.

\item Suppose that $\varphi: \fA \to \fB$ is a $\bw$-isomorphism. Then $\fB$ may be viewed as a right $\bw$-module over itself, and also as a left $\bw$-module over $\fA$ via the module action $a \cdot b = \varphi(a)b$ and inner product $\langle b|c \rangle_\fA = \varphi^{-1}(b c^*)$ for $a \in \fA$ and $b,c \in \fB$. It is straightforward to check that $\fB$ is an $\fA-\fB$ $\bw$-imprivitivity bimodule.

\item It is pretty direct to show that for any $n \in \N$, $\fB$ is $\bw$-Morita equivalent to $M_n(\fB)$ (the $\bw$-algebra of $n \times n$ matrices over $\fB$), via the $\bw$-module $C_n(\fB)$ (the $C^*$-module direct sum of $n$ copies of $\fB$). We show below (Corollary~\ref{fund Me examp}) that there is an infinite version of this, which is analogous to the facts that for any Hilbert space $\Hil$ of dimension $I$, any $C^*$-algebra $A$, and any $W^*$-algebra $M$, $A$ is strongly $C^*$-Morita equivalent to $\K(\Hil) \otimes A$ via $C_I(A)$, and $M$ is $W^*$-Morita equivalent to $B(\Hil) \overline{\otimes} M$ via $C_I^w(M)$. (Here, $C_I(A)$ denotes the $C^*$-module direct sum of $I$ copies of $A$, and $C_I^w(M)$ denotes the $W^*$-module direct sum of $I$ copies of $M$. 
See \cite[8.1.9 and 8.5.15]{BLM} for more details on these constructions.)
\end{enumerate}
\end{example}

We now show that condition (3) in Definition~\ref{Me} may be replaced with the condition that $\fA$ canonically embeds as a $\bw$-subalgebra of $\B_\fB(\fX)$.

\begin{lemma} \label{Me repn lem}
If $\fX$ is a $\sigma$-full left $\bw$-module over $\fA$ and a right $\bw$-module over $\fB$ such that ${}_\fA \langle x|y \rangle z = x \langle y|z \rangle_\fB$ for all $x,y,z \in \fX$, then there is a canonical isometric $*$-homomorphism $\lambda: \fA \hookrightarrow \B_\fB(\fX)$.
\end{lemma}

\proof
It is a well-known bit of basic $C^*$-module theory (cf.\ the first few lines in the proof of \cite[Lemma~8.1.15]{BLM}) that the inner product condition on $\fX$ implies that the canonical map $\lambda: \fA \to B(\fX)$, defined $\lambda(a)(x)=ax$ for $a \in \fA$ and $x \in \fX$, maps into $\B_\fB(\fX)$ and is a $*$-homomorphism.

We now use the assumption that $\fX$ is $\sigma$-full over $\fA$ to prove that $\lambda$ is injective, hence isometric. For this, suppose that $a \in \fA$ with $\lambda(a)=0$. Then $\lambda(a)(x)=ax=0$ for all $x \in \fX$, so $\langle ax|y \rangle_\fA = a \langle x|y \rangle_\fA=0$ for all $x,y \in \fX$. With $\mathscr{S}=\{j \in \fA: aj=0\}$, it is straightforward to check that $\langle \fX|\fX \rangle_\fA \subseteq \mathscr{S}$ and that $\mathscr{S}$ is $\sigma$-closed. Since $\langle \fX|\fX \rangle_\fA$ is $\sigma$-dense in $\fA$, $\mathscr{S}=\fA$. Hence $aa^*=0$, so $a=0$.
\endproof

\begin{lemma} \label{L nondeg gen lem}
If $S$ is any subset of $B(\Hil)$, then $[S \Hil] = [\Bw(S) \Hil]$. In particular, if $S$ is a $\sigma$-dense subset of a $\bw$-algebra $\fB \subseteq B(\Hil)$, then $[S \Hil]=\Hil$.
\end{lemma}

\proof
Obviously $[S \Hil] \subseteq [\Bw(S) \Hil]$. For the other inclusion, fix $\zeta \in \Hil$, and define $\mathscr{S}:=\{b \in \Bw(S) : b \zeta \in [S \Hil]\}$. Clearly $S \subseteq \mathscr{S}$, and by Mazur's lemma applied to $[S \Hil]$ (or a direct calculation), $\mathscr{S}$ is WOT sequentially closed. Thus $\Bw(S) = \mathscr S$, and the other inclusion follows. The last statement follows immediately (recalling our convention that concrete $\bw$-algebras are taken to be nondegenerately acting).
\endproof


Note that if $\fA \subseteq B(\Kil)$ is a $\bw$-algebra, then the copy of the multiplier algebra $M(\fA)$ in $B(\Kil)$ is WOT sequentially closed, hence is also a $\bw$-algebra. If $\fX$ is a left $\bw$-module over $\fA$, then $\fX$ is canonically a left $C^*$-module over $M(\fA)$ (see \cite[8.1.4~(4)]{BLM}), and it is easy to check that $\fX$ is a left $\bw$-module over $M(\fA)$ in this case.

\begin{lemma} \label{mult sigalg}
Let $\fA \subseteq B(\Kil)$ be a $\bw$-algebra, and let $(\xi_n)$ be a sequence in the copy of the multiplier algebra $M(\fA)$ in $B(\Kil)$. Then $(\xi_n)$ is WOT-convergent if and only if $(\xi_n a)$ is $\sigma$-convergent in $\fA$ for all $a \in \fA$. For $\xi \in M(\fA)$, we have $\xi_n \xrightarrow{WOT} \xi$ if and only if $\xi_n a \ssto \xi a$ for all $a \in \fA$. Hence $M(\fA)$ has a unique structure as a $\bw$-algebra in which $\fA$ is a $\bw$-subalgebra.
\end{lemma}

\proof
The forward direction of each of the first two claims is a direct consequence of separate WOT-continuity of the product in $B(\Kil)$. The converse of the first claim follows from a standard argument using Lemma~\ref{simplem} and the fact that vectors of the form $a \zeta$ for $a \in \fA$ and $\zeta \in \Kil$ are total in $\Kil$.

For the backward direction of the second claim, suppose $(\xi_n)$, $\xi \in M(\fA)$ such that $\xi_n a \ssto \xi a$ for all $a \in \fA$. By the first claim, there is a $\xi' \in M(\fA)$ such that $\xi_n \xrightarrow{WOT} \xi'$. Then $\xi a = \xi' a$ for all $a \in \fA$, which implies that $\xi = \xi'$ since $\fA$ is an essential ideal in $M(\fA)$. Hence $\xi_n \xrightarrow{WOT} \xi$.

For the final statement, suppose $\mathscr S$ is a $\sigma$-convergence system on $M(\fA)$ such that $(M(\fA), \mathscr S)$ is a $\bw$-algebra containing $\fA$ as a $\bw$-subalgebra. Let $\pi:(M(\fA), \mathscr S) \hookrightarrow B(\Hil)$ be a faithful $\bw$-representation. Then $\pi|_\fA$ is a faithful $\bw$-representation of $\fA$ since every nondegenerate representation of $M(\fA)$ restricts to a nondegenerate representation of $\fA$. (Indeed, to see the latter claim, suppose $\varphi:M(\fA) \to B(\Hil)$ is nondegenerate and $(e_\lambda)$ is a cai for $\fA$. Then $\varphi(e_\lambda)\varphi(\eta)\zeta \to \varphi(\eta)\zeta$ for all $\eta \in M(\fA)$ and $\zeta \in \Hil$. By Cohen's factorization theorem, $\{ \varphi(\eta) \zeta : \eta \in M(\fA), \, \zeta \in \Hil\} = \Hil$, so the claim follows.) Then $\xi_n \xrightarrow{\mathscr S} \xi$ in $M(\fA)$ if and only if $\pi(\xi_n) \xrightarrow{WOT} \pi(\xi)$ in $B(\Hil)$, which by what we proved above is equivalent to: $\pi(\xi_n a) \xrightarrow{WOT} \pi(\xi a)$ in $B(\Hil)$ for all $a \in \fA$. Since $\pi$ restricts to a faithful $\bw$-representation on $\fA$, the latter occurs if and only if $\xi_n a \ssto \xi a$ in $\fA$ for all $a \in \fA$, which by the above again is equivalent to $\xi_n \xrightarrow{WOT} \xi$ in $B(\Kil)$. So every such $\sigma$-convergence system $\mathscr S$ on $M(\fA)$ coincides with the one induced by the faithful $\bw$-representation $\fA \subseteq B(\Kil)$.
\endproof

The next lemma is a slight generalization of part of the previous lemma.

\begin{lemma} \label{L WOTconv}
Let $\fX$ be a $\sigma$-full left $\bw$-module over a $\bw$-algebra $\fA$, give $M(\fA)$ the $\bw$-algebra structure as in Lemma~\ref{mult sigalg}, and let $(\xi_n)$ be a sequence in $M(\fA)$. Then $(\xi_n)$ is $\sigma$-convergent in $M(\fA)$ if and only if $(\xi_n x)$ is ${}_\fA \sigma$-convergent for all $x \in \fX$. If $\xi \in M(\fA)$, $\xi_n \ssto \xi$ in $M(\fA)$ if and only if $\xi_n x \lmto{\fA} \xi x$ for all $x \in \fX$.
\end{lemma}

\proof
Suppose $(\xi_n)$ is $\sigma$-convergent in $M(\fA)$. Since $\xi_n \langle x|y \rangle = \langle \xi_n x|y \rangle$ for $x,y \in \fX$, the sequence $(\langle \xi_n x | y \rangle)$ is $\sigma$-convergent in $\fA$ for all $x, y \in \fX$. By Definition~\ref{ss-module def}, $(\xi_n x)$ is ${}_\fA \sigma$-convergent for all $x \in \fX$.

Conversely, suppose $(\xi_n x)$ is ${}_\fA \sigma$-convergent for all $x \in \fX$, and fix a faithful $\bw$-representation $\fA \subseteq M(\fA) \subseteq B(\Kil)$. A short calculation shows that $(\langle \xi_n j \zeta, \eta \rangle)$ converges for all $j \in \langle \fX|\fX \rangle$ and $\zeta, \eta \in \Kil$. The desired result follows by combining Lemma~\ref{L nondeg gen lem} and Lemma~\ref{simplem}.

The final statement follows similarly.
\endproof

Recall from \cite[Proposition~3.6]{Bear} that if $\fX$ is a right $\bw$-module over $\fB$, then $\B_\fB(\fX)$ is canonically a $\bw$-algebra in which $T_n \ssto T$ if and only if $T_n(x) \rmto{\fB} T(x)$ for all $x \in \fX$. It is this $\bw$-module structure on $\B_\fB(\fX)$ to which the results below refer.

\begin{prop} \label{Me equiv lr}
If $\fX$ is a $\bw$-Morita equivalence bimodule between $\fA$ and $\fB$, then the canonical isometric $*$-homomorphism $\lambda: \fA \hookrightarrow \B_\fB(\fX)$ from Lemma~\ref{Me repn lem} is a $\bw$-embedding.
\end{prop}

\proof
We first show that $\lambda$ is $\sigma$-continuous. Suppose that $a_n \ssto a$ in $\fA$. By Lemma~\ref{L WOTconv}, we have $\lambda(a_n)(x) = a_n x \lmto{\fA} ax = \lambda(a)(x)$ for all $x \in \fX$, so that $\lambda(a_n)(x) \rmto{\fB} \lambda(a)(x)$ for all $x \in \fX$ by Definition~\ref{Me}~(3). Hence $\lambda(a_n) \ssto \lambda(a)$ by the paragraph above the statement of the current proposition.

To see that $\lambda(\fA)$ is $\sigma$-closed, let $(a_n)$ be a sequence in $\fA$ such that $\lambda(a_n) \ssto T$ in $\B_\fB(\fX)$. Then $a_n x = \lambda(a_n)(x) \rmto{\fB} T(x)$, so $a_n x \lmto{\fA} T(x)$ for all $x \in \fX$. By Lemma~\ref{L WOTconv}, $(a_n)$ is $\sigma$-convergent to some $a \in \fA$. Hence $\lambda(a_n) \ssto \lambda(a)$, and thus $\lambda(\fA)$ is $\sigma$-closed.

To see that $\lambda^{-1}$ is $\sigma$-continuous, replace the $T$ in the previous paragraph with $\lambda(a')$. The argument there shows that there is some $a \in \fA$ with $a_n \ssto a$. So $\lambda(a')=\lambda(a)$, and since $\lambda$ is injective, $a'=a$ and $a_n \ssto a'$.
\endproof

\begin{prop} \label{Me equiv rl}
If $\fX$ is an $\fA-\fB$ bimodule satisfying (1) and (2) in Definition~\ref{Me} and such that the canonical isometric $*$-homomorphism $\lambda: \fA \hookrightarrow \B_\fB(\fX)$ from Lemma~\ref{Me repn lem} is a $\bw$-embedding, then $\fX$ is a $\bw$-imprivitivity bimodule.
\end{prop}

\proof
We need only to check (3) in Definition~\ref{Me}. Suppose that $x_n \lmto{\fA} x$ in $\fX$. Then for any $y, z \in \fX$,
\[ y \langle x_n|z \rangle_\fB = {}_\fA \langle y|x_n \rangle z = \lambda({}_\fA \langle y|x_n \rangle)(z) \rmto{\fB} \lambda({}_\fA \langle y|x \rangle)(z) = {}_\fA \langle y|x \rangle z = y \langle x_n|z \rangle_\fB.\]
By the ``right version" of Lemma~\ref{L WOTconv}, $\langle x_n|z \rangle_\fB \ssto \langle x|z \rangle_\fB$ for all $z \in \fX$; i.e., $x_n \rmto{\fB} x$. The other direction of (3) in the definition of $\bw$-imprivitivity bimodules follows by symmetry.
\endproof

\begin{thm} \label{Me equiv}
Let $\fX$ be an $\fA-\fB$ bimodule satisfying (1) and (2) in Definition~\ref{Me}. Then $\fX$ is a $\bw$-imprivitivity bimodule if and only if the canonical isometric $*$-homomorphism $\lambda: \fA \hookrightarrow \B_\fB(\fX)$ from Lemma~\ref{Me repn lem} is a $\bw$-embedding.

In this case, $\fA \cong \Bw(\K_\fB(\fX))$ $\bw$-isomorphically.
\end{thm}

\proof
Proposition~\ref{Me equiv lr} does one direction of the first claim, and Proposition~\ref{Me equiv rl} does the other. For the final claim, note that for $x,y \in \fX$, $\lambda({}_\fA \langle x|y \rangle)$ coincides with $|x \rangle \langle y|$ in $\K_\fB(\fX)$. Since the spans of these two types of elements generate $\lambda(\fA)$ and $\Bw(\K_\fB(\fX))$ as $\bw$-algebras respectively, we have that $\lambda(\fA) = \Bw(\K_\fB(\fX))$. Hence the claim is proved since $\lambda$ is a $\bw$-embedding by Proposition~\ref{Me equiv lr}.
\endproof


\section{The $\Sigma^*$-module interior tensor product of $\Sigma^*$-imprivitivity bimodules}

Let $\fX$ be a right $\bw$-module over $\fB$, let $\fY$ be a right $\bw$-module over $\fC$, and suppose there is a $\sigma$-continuous $*$-homomorphism $\lambda: \fB \to \B_\fC(\fY)$ with $\sigma$-closed range (although the definition below works if $\lambda$ is just a $*$-homomorphism). Recall the $C^*$-module interior tensor product $\fX \otimes_\fB \fY$ discussed in \ref{int tp constr} above---this is a right $C^*$-module over $\fC$. In direct analogy with the $W^*$-module case (see \cite[Section~3]{Ble97b}), we define the \emph{$\bw$-module interior tensor product} $\fX \otimes^{\bw}_\fB \fY$ to be the $\bw$-module completion (see \cite[3.11--3.14]{Bear}) of $\fX \otimes_\fB \fY$. In symbols,
\[
\fX \otimes^{\bw}_\fB \fY := \Bw(\fX \otimes_\fB \fY),
\]
where the $\sigma$-closure here is taken in $B(\Hil, (\fX \otimes_\fB \fY) \otimes_\fC \Hil)$ for any faithful $\bw$-representation $\fC \hookrightarrow B(\Hil)$.


For the rest of this section, let $\fA$, $\fB$, and $\fC$ be $\bw$-algebras, let $\fC \hookrightarrow B(\Hil)$ be a faithful $\bw$-representation, let $\fX$ be a $\bw$-imprivitivity bimodule between $\fA$ and $\fB$, and let $\fY$ be a $\bw$-imprivitivity bimodule between $\fB$ and $\fC$. By Proposition~\ref{Me equiv lr}, there is a canonical $\bw$-embedding $\fB \hookrightarrow \B_\fC(\fY)$, so we may take the $\bw$-module interior tensor product $\fX \otimes^{\bw}_\fB \fY$, which makes sense at least as a right $\bw$-module over $\fC$. 
The reader should keep in mind the inclusions
\[
\fX \otimes_\fB \fY \subseteq \fX \otimes^{\bw}_\fB \fY = \Bw(\fX \otimes_\fB \fY) \subseteq B(\Hil, (\fX \otimes_\fB \fY) \otimes_\fC \Hil).
\]

We now prove that $\fX \otimes^{\bw}_\fB \fY$ can be equipped with the structure of a left $\bw$-module over $\fA$ so that it becomes a $\bw$-imprivitivity bimodule between $\fA$ and $\fC$. The strategy is to show: (1) there is a canonical faithful $\bw$-representation $\pi: \fA \hookrightarrow B((\fX \otimes_\fB \fY) \otimes_\fC \Hil)$, so that $B(\Hil, (\fX \otimes_\fB \fY) \otimes_\fC \Hil)$ has a canonical left $\fA$-action; (2) $\fX \otimes^{\bw}_\fB \fY$ is an $\fA$-submodule of $B(\Hil, (\fX \otimes_\fB \fY) \otimes_\fC \Hil)$; (3) $\xi \eta^* \in \pi(\fA)$ for all $\xi, \eta \in \fX \otimes^{\bw}_\fB \fY$; (4) the $\fA$-valued map $\langle \xi|\eta \rangle := \pi^{-1}(\xi \eta^*)$ makes $\fX \otimes^{\bw}_\fB \fY$ into a left $\bw$-module over $\fA$; (5) with this structure, $\fX \otimes^{\bw}_\fB \fY$ is a $\bw$-imprivitivity bimodule between $\fA$ and $\fC$.

\begin{lemma} \label{L bw full}
With $\mathcal{J}:= \langle \fX \otimes_\fB \fY | \fX \otimes_\fB \fY \rangle_\fA$ and $\mathcal{I} :=  \langle \fX \otimes_\fB \fY | \fX \otimes_\fB \fY \rangle_\fC$, we have $\Bw(\mathcal{J}) = \fA$ and $\Bw(\mathcal{I}) = \fC$.
\end{lemma}

\proof
We prove $\Bw(\mathcal{J}) = \fA$; the other claim is similar. For $x,x' \in \fX$, and $y,y' \in \fY$, recall the formula $\langle x \otimes y|x' \otimes y' \rangle_\fA = \langle x \langle y|y' \rangle_\fB|x' \rangle_\fA$. It follows easily from this that $\mathcal{J} = \langle \fX \langle \fY|\fY \rangle_\fB | \fX \rangle_\fA$, where by the latter we mean  the span in $\fA$ of elements of the form $\langle x \langle y|y' \rangle_\fB x' \rangle_\fA$. With $x,x' \in \fX$ fixed, let $\mathscr{R}:=\{b \in \fB : \langle xb|x' \rangle_\fA \in \Bw(\mathcal{J})\}$. By the fact just mentioned, we easily see that $\langle \fY|\fY \rangle_\fB \subseteq \mathscr{R}$, and it follows from Lemma~\ref{L WOTconv} that $\mathscr{R}$ is $\sigma$-closed. Since $\fY$ is $\sigma$-full over $\fB$, we have $\mathscr{R} = \fB$. Since by Cohen's factorization theorem (\cite[A.6.2]{BLM}) we can write any $x \in \fX$ as $x=x'b$ for some $x' \in \fX$, $b \in \fB$, it follows that $\langle \fX | \fX \rangle_\fA \subseteq  \Bw(\mathcal{J})$. Thus $$\fA = \Bw(\langle \fX|\fX \rangle_\fA) \subseteq \Bw(\mathcal{J}) \subseteq \fA,$$ hence $\Bw(\mathcal{J}) = \fA$.
\endproof

\begin{lemma} \label{can emb}
There is a canonical faithful $\bw$-representation $\pi: \fA \hookrightarrow B((\fX \otimes_\fB \fY) \otimes_\fC \Hil)$.
\end{lemma}

\proof
By basic $C^*$-module theory, $\fX \otimes_\fB \fY$ is a left $C^*$-module over $\fA$ and right $C^*$-module over $\fC$, and a straightforward calculation (first checking on simple tensors) shows that $\langle \xi | \zeta \rangle_\fA \eta = \xi \langle \zeta|\eta \rangle_\fC$ for all $\xi, \zeta, \eta \in \fX \otimes_\fB \fY$. It follows as in the first few lines of Lemma~\ref{Me repn lem} there is a canonical $*$-homomorphism $\pi: \fA \to \B_\fC(\fX \otimes_\fB \fY) \subseteq B((\fX \otimes_\fB \fY) \otimes_\fC \Hil)$ To see that $\pi$ is injective, hence isometric, suppose that $a \in \fA$ with $\pi(a)=0$. Then
\[
a \langle x \otimes y |x' \otimes y' \rangle_\fA = \langle a(x \otimes y)|x' \otimes y' \rangle_\fA = \langle ax \otimes y|x' \otimes y' \rangle_\fA = 0
\]
for all $x,x' \in \fX$ and $y,y' \in \fY$. By taking linear combinations, we have $a j=0$ for all $j \in \mathcal{J} := \langle \fX \otimes_\fB \fY| \fX \otimes_\fB \fY \rangle_\fA$. The set $\mathscr{T}=\{ b \in \fA : ab = 0 \}$ is evidently $\sigma$-closed and contains $\mathcal{J}$, so $\mathscr{T} = \fA$ since $\Bw(\mathcal{J}) =\fA$ by Lemma~\ref{L bw full}. Hence $aa^*=0$, so $a=0$.

To see that $\pi$ is $\sigma$-continuous, suppose $a_n \ssto a$ in $\fA$. By Lemma~\ref{L WOTconv}, $a_n x \rmto{\fB} ax$ for all $x \in \fX$, hence by Lemma~\ref{L WOTconv} again, \[ \langle a_n x|x' \rangle_\fB y' \rmto{\fC} \langle a x|x' \rangle_\fB y' \] for all $x,x' \in \fX$ and $y' \in \fY$. A couple of easy calculations then give that for any $x,x' \in \fX$, $y, y' \in \fY$, and $k,k' \in \Hil$,
\begin{equation*}
\begin{split}
\langle \pi(a_n)((x \otimes y) \otimes k), (x' \otimes y') \otimes k' \rangle = \langle k, \langle y| \langle a_n x| x' \rangle_\fB y' \rangle_\fC k' \rangle  \\ \to \langle k, \langle y| \langle a x| x' \rangle_\fB y' \rangle_\fC k' \rangle = \langle \pi(a)((x \otimes y) \otimes k), (x' \otimes y') \otimes k' \rangle.
\end{split}
\end{equation*}
It is easy to check that elements of the form $(x \otimes y) \otimes k$ are total in $(\fX \otimes_\fB \fY) \otimes_\fC \Hil$, so by Lemma~\ref{simplem} we get $\pi(a_n) \xrightarrow{WOT} \pi(a)$.

We omit these, but a couple of routine arguments using Lemma~\ref{L nondeg gen lem} and Lemma~\ref{simplem} show that $\pi$ has WOT sequentially closed range and that $\pi^{-1}$ is $\sigma$-continuous.
\endproof

\begin{lemma}
Let $\pi$ be as in the previous lemma, and view $B(\Hil,(\fX \otimes_\fB \fY)\otimes_\fC \Hil)$ as a left $\fA$-module via the action $a \cdot T = \pi(a)T$. Then $\bwtp$ is an $\fA$-submodule of $B(\Hil,(\fX \otimes_\fB \fY)\otimes_\fC \Hil)$.
\end{lemma}

\proof
Fix $a \in \fA$. Then, as mentioned at the beginning of the proof of the previous lemma, $\pi(a)$ lies in the canonical copy of $\B_\fC(\fX \otimes_\fB \fY)$ in $B((\fX \otimes_\fB \fY) \otimes_\fB \Hil)$, so it follows that for any $\eta$ in $\fX \otimes_\fB \fY$, $\pi(a) \eta$ remains in $\fX \otimes_\fB \fY$. We now employ our usual trick: setting $\mathscr{R}=\{ \eta \in \bwtp : \pi(a) \eta \in \bwtp \}$, we have $\fX \otimes_\fB \fY \subseteq \mathscr{R}$, and the latter is WOT sequentially closed since $\bwtp$ is. So $\mathscr{R} = \bwtp$.
\endproof

\begin{lemma} \label{L ipvals}
Let $\pi$ be as in the previous lemma. If $\xi, \eta \in \bwtp \subseteq B(\Hil, (\fX \otimes_\fB \fY) \otimes_\fC \Hil)$, then $\xi \eta^* \in \pi(\fA)$.
\end{lemma}

\proof
First suppose that $\xi, \eta \in \fX \otimes_\fB \fY$, which as we mentioned at the beginning of the proof of Lemma~\ref{can emb} is a left $C^*$-module over $\fA$ and right $C^*$-module over $\fC$ satisfying the inner product relation $\langle \xi|\zeta \rangle_\fA \eta = \xi \langle \zeta|\eta \rangle_\fC$. Then for any $\zeta, \zeta' \in \fX \otimes_\fB \fY$ and $k,k' \in \Hil$, we have
\begin{align*}
\langle \xi \eta^*(\zeta \otimes k), \zeta' \otimes k' \rangle
	&= \langle \eta^* (\zeta \otimes k), \xi^*(\zeta' \otimes k') \rangle \\
	&= \langle \langle \eta|\zeta \rangle_\fC k, \langle \xi|\zeta' \rangle_\fC  k' \rangle	\\
	&= \langle \xi \otimes \langle \eta|\zeta \rangle_\fC k, \zeta' \otimes k' \rangle \\
	&= \langle \xi \langle \eta|\zeta \rangle_\fC \otimes k, \zeta' \otimes k' \rangle \\
	&= \langle \langle \xi|\eta \rangle_\fA \zeta \otimes k,\zeta' \otimes k' \rangle \\
	&= \langle \pi(\langle \xi|\eta \rangle_\fA)(\zeta \otimes k), \zeta' \otimes k' \rangle.
\end{align*}
Since the simple tensors are total, we conclude that $\xi \eta^* = \pi(\langle \xi|\eta \rangle_\fA) \in \pi(\fA)$.

Keep $\xi \in \fX \otimes_\fB \fY$, and define $\mathscr{R}:=\{\eta \in \bwtp: \xi \eta^* \in \pi(\fA)\}$ (recalling again that we are viewing $\bwtp$ as the WOT sequential closure of $\fX \otimes_\fB \fY$ in $B(\Hil, (\fX \otimes_\fB \fY) \otimes_\fC \Hil)$). By what we just proved, $\fX \otimes_\fB \fY \subseteq \mathscr{R}$, and since $\pi(\fA)$ is WOT sequentially closed by Lemma~\ref{can emb}, it follows that $\mathscr{R}$ is WOT sequentially closed, so $\mathscr{R}=\bwtp$.

A similar argument using this shows the full claim.
\endproof

\begin{thm} \label{P Metrans}
With $\fA$-module action $a \cdot \xi = \pi(a) \xi$ and $\fA$-valued inner product defined $\langle \xi|\eta \rangle_\fA := \pi^{-1}(\xi \eta^*)$, $\fX \otimes_\fB^{\bw} \fY$ is a $\bw$-imprivitivity bimodule between $\fA$ and $\fC$.
\end{thm}

\proof
We have already mentioned that the work in \cite{Bear} already gives that $\bwtp$ is a right $\bw$-module over $\fC$.

We now describe why $\bwtp$ is a left $C^*$-module over $\fA$. Most of the axioms defining a $C^*$-module tensor product are easily checked for $\langle \cdot|\cdot \rangle_\fA$ using the fact that $\pi^{-1}$ is a homomorphism. Completeness of $\bwtp$ in the norm induced by $\langle \cdot|\cdot \rangle_\fA$ follows since $\bwtp$ is complete in the norm it inherits from $B((\fX \otimes_\fB \fY)\otimes_\fC \Hil)$, and since these two norms coincide:
\[ \|\xi\|^2_{B((\fX \otimes_\fB \fY)\otimes_\fC \Hil)} = \| \xi \xi^* \|_{B((\fX \otimes_\fB \fY)\otimes_\fC \Hil)} = \|\pi^{-1}(\xi \xi^*)\|_\fA = \|\langle \xi|\xi \rangle_\fA\| \]
for $\xi \in \bwtp$, using the fact that $\pi^{-1}$ is isometric.

To see that $\bwtp$ is a left $\bw$-module over $\fA$, suppose that $(\xi_n)$ is a sequence in $\bwtp$ such that $\langle \xi_n|\eta \rangle_\fA$ is $\sigma$-convergent for all $\eta \in \bwtp$. Then $\xi_n \eta^* = \pi(\langle \xi_n|\eta \rangle)$ is WOT-convergent in $B((\fX \otimes_\fB \fY) \otimes_\fC \Hil)$ since $\pi$ is a $\bw$-representation. By Lemma~\ref{L bw full}, $\Bw(\langle \fX \otimes_\fB \fY|\fX \otimes_\fB \fY \rangle_\fC \Hil) = \fC$. So by Lemma~\ref{L nondeg gen lem}, $[(\fX \otimes_\fB \fY)^*((\fX \otimes_\fB \fY) \otimes_\fC \Hil)] = [\langle \fX \otimes_\fB \fY|\fX \otimes_\fB \fY \rangle_\fC \Hil] = \Hil$, which implies that elements of the form $\eta^* h$, with $\eta \in \fX \otimes_\fB \fY$ and $h \in (\fX \otimes_\fB \fY) \otimes_\fC \Hil$, are total in $\Hil$. Since $\langle \xi_n \eta^* h,k \rangle$ converges for all $\eta \in \fX \otimes_\fB \fY$ and $h,k \in (\fX \otimes_\fB \fY) \otimes_\fC \Hil$, we may conclude by Lemma~\ref{simplem} that $(\xi_n)$ WOT-converges to some $\xi \in B(\Hil,(\fX \otimes_\fB \fY) \otimes_\fC \Hil)$. This $\xi$ is in $\bwtp$ since the latter is WOT sequentially closed, and we have $$\langle \xi_n|\eta \rangle_\fA = \pi^{-1}(\xi_n \eta^*) \ssto \pi^{-1}(\xi \eta^*) = \langle \xi|\eta \rangle_\fA \text{ for all } \eta \in \bwtp$$ since $\pi^{-1}$ is $\sigma$-continuous. Hence $\bwtp$ is a $\bw$-module.

We now check condition (3) in Definition~\ref{Me}. We showed in the previous paragraph that if $\xi_n \lmto{\fA} \xi$, then $\xi_n \xrightarrow{WOT} \xi$ in $B(\Hil,(\fX \otimes_\fB \fY) \otimes_\fC \Hil)$, but the latter is equivalent to $\xi_n \rmto{\fC} \xi$. Conversely, if $\xi_n \rmto{\fC} \xi$, then $\xi_n \eta^* \xrightarrow{WOT} \xi \eta^*$ in $B((\fX \otimes_\fB \fY) \otimes_\fC \Hil)$ for all $\eta \in B(\Hil,(\fX \otimes_\fB \fY) \otimes_\fC \Hil)$. In particular, since $\pi^{-1}$ is $\sigma$-continuous, $\langle \xi_n|\eta \rangle_\fA = \pi^{-1}(\xi_n \eta^*) \ssto \pi^{-1}(\xi \eta^*) = \langle \xi|\eta \rangle_\fA$ for all $\eta \in \bwtp$. Hence $\xi_n \lmto{\fA} \xi$.

Lemma~\ref{L bw full} shows that $\bwtp$ is $\bw$-full over both $\fA$ and $\fB$ (the former since the $\fA$-valued inner product on $\bwtp$ extends the one on $\fX \otimes_\fB \fY$ by construction), so it only remains to check the inner product formula in Definition~\ref{Me}~(3). This is straightforward from the definitions of the $\fA$-module action and inner product:
\[ \langle \xi|\eta \rangle_\fA \zeta = \pi^{-1}(\xi \eta^*) \cdot \zeta = (\xi \eta^*) \zeta = \xi(\eta^* \zeta) = \xi \langle \eta|\zeta \rangle_\fC\]
for $\xi,\eta, \zeta \in \bwtp$.
\endproof

\begin{cor} \label{T er}
Strong $\bw$-Morita equivalence is an equivalence relation strictly coarser than $\bw$-isomorphism.
\end{cor}

\proof
Reflexivity follows from (2) in Example~\ref{ssMe examps}.

Symmetry follows just as in $C^*$-Morita equivalence---if $\fX$ is an $\fA-\fB$ $\bw$-imprivitivity bimodule, then the adjoint $C^*$-module $\overline{\fX}$ (see \cite[8.1.1]{BLM}) is easily checked to be a $\fB-\fA$ $\bw$-imprivitivity bimodule.

Transitivity was proved in Theorem~\ref{P Metrans}.

The ``coarser" claim was shown in Example~\ref{ssMe examps}~(2), and the ``strictly" part may be seen from Example~\ref{ssMe examps}~(3).
\endproof

\section{Weak $\bw$-Morita equivalence}

In his original paper \cite{Rie} on the subject, Rieffel proved that two $W^*$-algebras are $W^*$-Morita equivalent if and only if their categories of normal Hilbert space representations are equivalent (the latter is actually his definition of Morita equivalence, see \cite[Definition~7.4]{Rie}; the result is \cite[Theorem~7.9]{Rie}). The corresponding statement for $C^*$-algebras is however not true, and so there are two different notions of $C^*$-algebraic Morita equivalence---the ``strong" version we have already defined, and a ``weak" version described below. After describing weak $C^*$-Morita equivalence, we define and study the $\bw$-version, which we call ``weak $\bw$-Morita equivalence". This turns out (as in the $C^*$-case and differing from the $W^*$-case) to be strictly coarser than strong $\bw$-Morita equivalence.

Let $A$ be a $C^*$-algebra. A \emph{Hilbert $A$-module} is a Hilbert space $H$ carrying an $A$-module structure via a $*$-homomorphism $A \to B(H)$. Let ${}_A HMOD$ be the category of Hilbert $A$-modules with morphisms the bounded $A$-module maps. For $C^*$-algebras $A$ and $B$, a \emph{$*$-functor} from ${}_A HMOD$ to ${}_B HMOD$ is a linear functor $F: {}_A HMOD \to {}_B HMOD$ such that $F(T^*)=F(T)^*$ whenever $T : H \to K$ is a morphism of Hilbert $A$-modules. Recall that two categories $\mathscr C$ and $\mathscr D$ are \emph{equivalent} if there are functors $F: \mathscr C \to \mathscr D$ and $G : \mathscr D \to \mathscr C$ with natural isomorphisms $FG \cong \mathrm{id}_\mathscr C$ and $GF \cong \mathrm{id}_\mathscr D$ (see \cite[IV.4]{SML}).

\begin{defn}
Two $C^*$-algebras $A$ and $B$ are \emph{weakly $C^*$-Morita equivalent} if ${}_A HMOD$ and ${}_B HMOD$ are equivalent via $*$-functors.
\end{defn}

\begin{note}
Weak $C^*$-Morita equivalence is indeed no stronger than strong $C^*$-Morita equivalence  (\cite[Theorem~6.23]{Rie2}). As mentioned above though, two weakly $C^*$-Morita equivalent $C^*$-algebras need not be strongly $C^*$-Morita equivalent. For example, $C([0,1])$ and $C(\mathbb T)$ are weakly but not strongly $C^*$-Morita equivalent. (Indeed, two commutative $C^*$-algebras are strongly $C^*$-Morita equivalent if and only if they are $*$-isomorphic (\cite[Corollary~6.27]{Rie2}), but any two separable commutative $C^*$-algebras whose spectra have the same cardinality are weakly $C^*$-Morita equivalent (\cite[Proposition~8.18]{Rie}).) 
See also Beer's paper \cite{Beer} for further interesting results about weak and strong $C^*$-Morita equivalence. A word of warning on the terminology: in Beer's and Rieffel's works, weak $C^*$-Morita equivalence is called simply ``Morita equivalence".
\end{note}

We now investigate the $\bw$-analogue of weak $C^*$-Morita equivalence, showing among other things that weak $\bw$-Morita equivalence is an equivalence relation strictly coarser than strong $\bw$-Morita equivalence.

Let $\fA$ be a $\bw$-algebra. If $H$ is a Hilbert space carrying an $\fA$-module structure via a $\bw$-representation $\fA \to B(H)$, say that $H$ is a \emph{$\bw$-Hilbert $\fA$-module}. Let ${}_\fA \bw HMOD$ be the category of $\bw$-Hilbert $A$-modules with morphisms the bounded $\fA$-module maps.

\begin{defn}
Two $\bw$-algebras $\fA$ and $\fB$ are weakly $\bw$-Morita equivalent if ${}_\fA \bw HMOD$ and ${}_\fB \bw HMOD$ are equivalent via $*$-functors.
\end{defn}

If $A$ is a $C^*$-algebra, its \emph{Davies-Baire envelope}, denoted $\Sigma(A)$, is the WOT sequential closure of $A$ in its universal representation. (Note that $\Sigma(A)$ may be identified with the weak* sequential closure of $A$ in $A^{**}$.)

\begin{prop} \label{weak me for envs}
If $A$ and $B$ are $C^*$-algebras, the following are equivalent:
\begin{enumerate}
\item $A$ and $B$ are weakly $C^*$-Morita equivalent;
\item $\Sigma(A)$ and $\Sigma(B)$ are weakly $\bw$-Morita equivalent;
\item $A^{**}$ and $B^{**}$ are $W^*$-Morita equivalent.
\end{enumerate}
\end{prop}

\proof
Equivalence of the first two follows directly from the correspondences between representations of a $C^*$-algebra and $\bw$-representations of its Davies-Baire envelope (see \cite{Dav68}, Theorem~3.1 and its proof). Equivalence of the first and last follows similarly (\cite[Proposition~8.18]{Rie})
\endproof

\begin{lemma} \label{sshmod lem}
Suppose $\fA$ and $\fB$ are $\bw$-algebras and $\fX$ is an $\fA-\fB$ $\bw$-imprivitivity bimodule. If $\Kil$ is a $\bw$-Hilbert $\fA$-module, then $\overline{\fX} \otimes_\fA \Kil$ is a $\bw$-Hilbert $\fB$-module. Similarly, if $\Hil$ is a $\bw$-Hilbert $\fB$-module, then $\fX \otimes_\fB \Hil$ is a $\bw$-Hilbert $\fA$-module.
\end{lemma}

\proof
Routine calculations show that if $\Kil$ is a $\bw$-Hilbert $\fA$-module, then there is a $*$-homomorphism $\pi: \fB \to B(\overline{\fX} \otimes_\fA \Kil)$ determined by the formula
\[ \pi(b)(\overline{x} \otimes \zeta) = \overline{x b^*} \otimes \zeta \] on simple tensors. To see that $\pi$ is $\sigma$-continuous, suppose $b_n \ssto b$ in $\fB$. Then $(\pi(b_n))$ is a bounded sequence, and for any $x, y \in \fX$ and $\zeta, \eta \in \Kil$,
\[ \langle \pi(b_n)(\overline x \otimes \zeta), \overline y \otimes \eta \rangle = 
\langle \zeta, {}_\fA \langle x b_n |y \rangle \eta \rangle \to \langle \zeta, {}_\fA \langle xb|y \rangle \eta \rangle = \langle \pi(b) (\overline x \otimes \zeta), \overline y \otimes \eta \rangle \]
since $x b_n \lmto{\fA} xb$ by Lemma~\ref{L WOTconv} and since $\fX$ is a $\bw$-imprivitivity bimodule. Hence $\pi(b_n) \xrightarrow{WOT} \pi(b)$ by Lemma~\ref{simplem}.

The other claim follows similarly.
\endproof

\begin{prop}
If $\fA$ and $\fB$ are strongly $\bw$-Morita equivalent $\bw$-algebras, then they are weakly $\bw$-Morita equivalent.
\end{prop}

\proof
Let $\fX$ be an $\fA-\fB$ $\bw$-imprivitivity bimodule. Define a functor $$F: {}_\fA \bw HMOD \to {}_\fB \bw HMOD$$ by $F(\Kil) = \overline{\fX} \otimes_\fA \Kil$ for a $\bw$-Hilbert $\fA$-module $\Kil$ and $F(T) = \mathrm{id}_{\overline{\fX}} \otimes T$ for a bounded $\fA$-module map $T: \Kil \to \Kil'$ between $\bw$-Hilbert $\fA$-modules. (The existence of $\mathrm{id}_{\overline{\fX}} \otimes T$ follows from functoriality of the $C^*$-module interior tensor product---see \cite[8.2.12~(1)]{BLM}.) Similarly, define a functor $$G: {}_\fB \bw HMOD \to {}_\fA \bw HMOD$$ by $G(\Hil) = \fX \otimes_\fB \Hil$ and $G(S) = \mathrm{id_\fX} \otimes S$.

By Lemma~\ref{sshmod lem}, these functors do indeed map into the desired categories, and it is straightforward to check that they are $*$-functors. To see that $FG$ is naturally isomorphic to $\mathrm{id}_{{}_\fB \bw HMOD}$, take a $\bw$-Hilbert $\fB$-module $\Hil$. Then by standard facts about the $C^*$-module interior tensor product (or see \cite[8.2.19]{BLM}), we have canonical $\fB$-module isomorphisms
\begin{align*}
FG(\Hil) &= \overline{\fX} \otimes_\fA (\fX \otimes_\fB \Hil)\\ 
&\cong (\overline{\fX} \otimes_\fA \fX) \otimes_\fB \Hil\\ &\cong {}_\fA \K(\fX) \otimes_\fB \Hil\\
&\cong \langle \fX| \fX \rangle_\fB \otimes_\fB \Hil\\
&\cong [\langle \fX| \fX \rangle_\fB \Hil]\\
&= \Hil,
\end{align*}
where the last line holds by Lemma~\ref{L nondeg gen lem}.
(To see the fourth line, note that since $\fX$ is a $\langle \fX|\fX \rangle_\fA-\langle \fX|\fX \rangle_\fB$ $C^*$-imprivitivity bimodule, ${}_{\langle \fX|\fX \rangle_\fA} \K(\fX) = {}_\fA \K(\fX)$ with its canonical $C^*$-algebra structure is canonically $*$-isomorphic to $\langle \fX | \fX \rangle_\fB$, and it is easy to check that this $*$-isomorphism is a $\fB$-$\fB$ bimodule map.)

The induced transformation from $FG$ to $\mathrm{id}_{{}_\fB \bw HMOD}$ is easily check to be a natural isomorphism. That $GH$ is naturally isomorphic to $\mathrm{id}_{{}_\fA \bw HMOD}$ may be checked similarly.
\endproof

To show that weak $\bw$-Morita equivalence is strictly weaker than strong $\bw$-Morita equivalence, it suffices as in the $C^*$-case to look at commutative algebras, but one has to look a little deeper.

For a locally compact Hausdorff space $X$, denote by $\mathrm{Baire}(X)$ the space of complex-valued Baire functions on $X$ in the sense of \cite[6.2.10]{Ped NOW}. That is, $\mathrm{Baire}(X)$ is the bounded pointwise sequential closure of $C_0(X)$ inside the space of all bounded functions on $X$, so that $\mathrm{Baire}(X)$ is a concrete $\bw$-algebra in $B(\ell^2(X))$ (in which WOT-convergence of sequences coincides with pointwise convergence of bounded sequences of functions). If $X$ is second countable, $\mathrm{Baire}(X)$ coincides with the space of bounded Borel measurable functions on $X$. A short exercise shows that $X$ is $\sigma$-compact if and only if $\mathrm{Baire}(X)$ is unital.

The proof of the following is essentially the same as one that works in the $C^*$-case (see \cite{BLM}, note on top of pg.~352).

\begin{prop} \label{strong me for comm}
Commutative $\bw$-algebras are strongly $\bw$-Morita equivalent if and only if they are $\bw$-isomorphic.
\end{prop}

\proof
Let $\fA$ and $\fB$ be commutative $\bw$-algebras, and let $\fX$ be an $\fA-\fB$ $\bw$-imprivitivity bimodule. We will show that there is a $\bw$-isomorphism $\varphi:\fA \to \fB$ determined by the formula $ax=x \varphi(a)$ for $x \in \fX$, $a \in \fA$.

Consider $\fX$ as a left $\bw$-module over $M(\fA)$ and a right $\bw$-module over $M(\fB)$ (see discussion above Lemma~\ref{mult sigalg}). We claim that for any fixed $\eta \in M(\fA)$, there is a unique $\zeta \in M(\fB)$ such that $\eta x = x \zeta$ for all $x \in \fX$. A short calculation using commutativity of $M(\fA)$ shows that the map $x \mapsto \eta x$ is in ${}_\fA \mathbb{B}(\fX) = M({}_\fA \K(\fX)) \subseteq M(\Bw({}_\fA \K(\fX)))$, so the canonical identification of the latter with $M(\fB)$ from the left-module version of Theorem~\ref{Me equiv} gives existence. For uniqueness, suppose $\zeta_1, \zeta_2 \in M(\fB)$ satisfy $x \zeta_1=x \zeta_2$ for all $x \in \fX$. Fixing $b \in \fB$, we have $xb(\zeta_1-\zeta_2) = 0$ for all $x \in \fX$, so that $\langle y|x \rangle b(\zeta_1-\zeta_2)=0$ for all $y,x \in \fX$. Hence $\mathcal J:=\{c \in \fB : cb(\zeta_1-\zeta_2)=0\}$ is a WOT sequentially closed subset of $\fB$ containing $\langle \fX|\fX \rangle_\fB$, so $\mathcal J = \fB$ and $b(\zeta_1-\zeta_2)=0$. Thus $\fB(\zeta_1-\zeta_2)=0$, and since $\fB$ is an essential ideal in $M(\fB)$, $\zeta_1=\zeta_2$.

We have proved the existence of a map $\theta: M(\fA) \to M(\fB)$ determined by the formula $\eta x = x \theta(\eta)$ for $x \in \fX$, $\eta \in M(\fA)$. By symmetry, $\theta$ has an inverse, so is injective and surjective. It is direct to check that $\theta$ is in fact a $*$-isomorphism (part of this uses commutativity of $M(\fA)$ again). To see that $\theta$ is a $\bw$-isomorphism, suppose that $(\eta_n)$ is a sequence in $M(\fA)$ with $\eta_n \ssto \eta$. Lemma~\ref{L WOTconv} shows that $x \theta(\eta_n) = \eta_n x \lmto{\fA} \eta x = x \theta(\eta)$ for all $x \in \fX$. Hence $x \theta(\eta_n) \rmto{\fB} x \theta(\eta)$ for all $x \in \fX$ since $\fX$ is a $\bw$-imprivitivity bimodule. By the ``right version" of Lemma~\ref{L WOTconv} again, we have $\theta(\eta_n) \ssto \theta(\eta)$. By symmetry, $\theta^{-1}$ is also $\sigma$-continuous.

It remains to show that $\theta(\fA) = \fB$. (This part of the proof is essentially from the end of the proof of 8.6.5 in \cite{BLM}.) For $x,y \in \fX$, we have
\begin{align*}
{}_\fA \langle y|y \rangle^2 x &= y \langle y|y \rangle_\fB \langle y|x \rangle_\fB
= y \langle y|x \rangle_\fB \langle y|y \rangle_\fB \\
&= {}_\fA \langle y|y \rangle {}_\fA \langle x|y \rangle y
= {}_\fA \langle x|y \rangle {}_\fA \langle y|y \rangle y \\
&= x \langle y|y \rangle_\fB^2
\end{align*}
by Definition~\ref{Me}~(2) and commutativity of $\fA$ and $\fB$. Hence $$\theta({}_\fA \langle y|y \rangle)^2 = \theta({}_\fA \langle y|y \rangle^2) = \langle y|y \rangle_\fB^2,$$ so $\theta({}_\fA \langle y|y \rangle)=\langle y|y \rangle_\fB$. By polarization, we have $\theta({}_\fA \langle y|z \rangle)= \langle z|y \rangle_\fB$ for all $y,z \in \fX$, and it follows that $\theta({}_\fA \langle \fX|\fX \rangle) = \langle \fX|\fX \rangle_\fB$. Since $\theta$ is $\sigma$-continuous, $\theta(\fA) \subseteq \fB$. A symmetric argument shows that $\theta^{-1}(\fB) \subseteq \fA$. Hence $\theta(\fA)=\fB$ and $\varphi = \theta|_\fA$ is a $\bw$-isomorphism between $\fA$ and $\fB$.
\endproof

\begin{example}
Let $L$ denote the long line (see \cite[Section~24, Exercise~12]{Munkres}). By \cite[pg.~257]{Beer}, $C([0,1])$ and $C_0(L)$ are weakly $C^*$-Morita equivalent. Hence by Proposition~\ref{weak me for envs}, $\mathrm{Baire}([0,1])$ and $\mathrm{Baire}(L)$ are weakly $\bw$-Morita equivalent. On the other hand, $L$ is not $\sigma$-compact (an exercise in basic topology shows that $L$ is not even Lindel\"of), so $\mathrm{Baire}(L)$ is non-unital, and thus $\mathrm{Baire}([0,1])$ and $\mathrm{Baire}(L)$ are not $*$-isomorphic. Hence they are not strongly $\bw$-Morita equivalent by Proposition~\ref{strong me for comm}.
\end{example}

\section{The ``full corners" characterization}

\begin{defn}
Let $\fC$ be a $\bw$-algebra. Two $C^*$-subalgebras $\fA$ and $\fB$ in $\fC$ are called \emph{complementary $\sigma$-full corners of $\fC$} if there is a projection $p \in M(\fC)$ such that $p \fC p \cong \fA$ $\bw$-isomorphically, $p^\perp \fC p^\perp \cong \fB$ $\bw$-isomorphically, and $\Bw(\fC p \fC) = \fC = \Bw(\fC p^\perp \fC)$.
\end{defn}

\begin{thm} \label{T corners are me}
Two $\bw$-algebras $\fA$ and $\fB$ are strongly $\bw$-Morita equivalent if and only if they are $\bw$-isomorphic to complementary $\sigma$-full corners of a $\bw$-algebra $\fC$. If $p \fC p$ and $p^\perp \fC p^\perp$ are complementary $\sigma$-full corners in a $\bw$-algebra $\fC$, then $p \fC p^\perp$ is a $p \fC p-p^\perp \fC p^\perp$ $\bw$-imprivitivity bimodule.
\end{thm}

\proof
$(\implies)$ Let $\fX$ be an $\fA-\fB$ $\bw$-imprivitivity bimodule. Define $\blink \fX$ to be the space of $2 \times 2$ matrices with entries as follows: \[ \blink \fX := \left[ \begin{matrix} \Bw(\K_\fB(\fX)) & \fX \\ \overline{\fX} & \fB \end{matrix} \right]. \] By \cite[Proposition~3.9]{Bear}, $\blink \fX$ (with the canonical product and adjoint arising from formal $2 \times 2$ matrix multiplication), $\blink \fX$ is canonically a $\bw$-algebra. Setting $p:=\left[ \begin{matrix} 1 & 0 \\ 0 & 0 \end{matrix} \right]$, we have $\fB = p^\perp \blink \fX p^\perp$ and $\fA \cong \Bw(\K_\fB(\fX)) = p \blink \fX p$ $\bw$-isomorphically by the last statement in Theorem~\ref{Me equiv}. It is easy to check using Cohen's factorization theorem and $\sigma$-fullness of $\fX$ that $\Bw(\blink \fX p \blink \fX) = \blink \fX$ and $\Bw(\blink \fX p^\perp \blink \fX) = \blink \fX$.

$(\impliedby)$ Assume that $\fA$ and $\fB$ are $\bw$-isomorphic to complementary $\sigma$-full corners of a $\bw$-algebra $\fC$, and set $\fX := p \fC p^\perp$. We will show that $\fX$ is a $\bw$-imprivitivity bimodule between $p \fC p$ and $p^\perp \fC p^\perp$ (thus also proving the final statement). As observed in \cite[Theorem~3.10]{Bear}, for a faithful $\bw$-representation $\fC \hookrightarrow B(\Hil)$, $\fX$ is a left (resp.\ right) $\bw$-module over the $\bw$-algebra $p \fC p \subseteq B(p \Hil)$ (resp.\ $p^\perp \fC p^\perp \subseteq B(p^\perp \Hil)$). To see that $\fX$ is $\sigma$-full over these, note that $\Bw_{p \Hil}(\langle \fX|\fX \rangle_{p \fC p}) = \Bw_{p \Hil}(p \fC p^\perp \fC p) = p \Bw_{\Hil}(\fC p^\perp \fC)p = p \fC p$ and similarly $\Bw_{p^\perp \Hil}(\langle \fX|\fX \rangle_{p^\perp \fC p^\perp}) = p^\perp \fC p^\perp$ (using the easy fact that $\Bw_{p \Hil}(p A p) = p \Bw_\Hil(A) p$ for any $C^*$-algebra $A \subseteq B(\Hil)$). So we have checked condition (1) in Definition~\ref{Me}. Condition (2) is obvious, so it remains to check condition (3); i.e., that for $(x_n),x \in \fX$, $x_n \lmto{p \fC p} x$ if and only if $x_n \rmto{p^\perp \fC p^\perp} x$. Note first that $[p^\perp \fC p \Hil] = p^\perp \Hil$ and $[p \fC p^\perp \Hil] = p \Hil$. Indeed, since $\Bw(p^\perp \fC p \fC) = p^\perp \fC$, we have
\[ [p^\perp \fC p \Hil] = [p^\perp \fC p \fC \Hil ] = [ p^\perp \fC \Hil] = p^\perp \Hil \]
by Lemma~\ref{L nondeg gen lem}, and the other equation is proved similarly. The desired result then follows readily from straightforward calculations and Lemma~\ref{simplem}.
\endproof

\begin{note}[Rieffel subequivalence for $\bw$-algebras]
We will not provide the details, but by following \cite[8.2.24]{BLM}, replacing norm-closures with $\sigma$-closures, one may prove that if $\fA$ and $\fB$ are strongly $\bw$-Morita equivalent via a $\bw$-imprivitivity bimodule $\fX$, then there are lattice isomorphisms between the following: (1) the $\sigma$-closed two-sided ideals of $\fA$, (2) the $\sigma$-closed two-sided ideals of $\fB$, (3) the ${}_\fA \sigma$-closed (equivalently, $\sigma_\fB$-closed) $\fA-\fB$ submodules of $\fX$, and (4) the $\sigma$-closed two-sided ideals of $\blink \fX$.
\end{note}

\begin{note}[The TRO picture]
Recall that a \emph{ternary ring of operators} (\emph{TRO} for short) is a closed subspace $X \subseteq B(\Hil,\Kil)$ for Hilbert spaces $\Hil$, $\Kil$, such that $xy^*z \in X$ for all $x,y,z \in X$. Just as Theorem~\ref{T corners are me} is a ``corners picture" of strong $\bw$-Morita equivalence, there is also a ``TRO picture" of the same. Namely, two $\bw$-algebra $\fA$ and $\fB$ are strongly $\bw$-Morita equivalent if and only if there is a WOT sequentially closed TRO $\fX \subseteq B(\Hil,\Kil)$ for some $\Hil, \Kil$, such that $\Bw(\fX \fX^*) \cong \fA$ and $\Bw(\fX^* \fX) \cong \fB$ $\bw$-isomorphically. (This is similar to the ``corner picture" and ``TRO picture" of $\bw$-modules presented in \cite[Theorem~3.10]{Bear}.)

We will omit most of the details of the proof since they are direct and much in the same line as other proofs already given here. For the forward direction, one fixes a faithful $\bw$-representation $\fB \hookrightarrow B(\Hil)$ and takes $\fX$ to be the image of a $\bw$-imprivitivity bimodule in $B(\Hil, \fX \otimes_\fB \Hil)$. The converse follows from the fact that if $\fX$ is a WOT sequentially closed TRO, then $\fX$ is a $\bw$-imprivitivity bimodule between $\Bw(\fX \fX^*)$ and $\Bw(\fX^* \fX)$.
\end{note}

Let $\fB \subseteq B(\Hil)$ be a $\bw$-algebra, let $I$ be a cardinal number, and denote by $\K_I$ the $C^*$-algebra of compact operators on a Hilbert space of dimension $I$. Recalling that the $C^*$-algebra tensor product $\K_I \otimes \fB$ may be viewed as a certain space of infinite matrices over $\fB$  (see e.g.\ \cite[1.2.26 and 1.5.2]{BLM}), we prefer to denote this space as $\K_I(\fB)$, and we view it canonically as a concrete $C^*$-algebra in $B(\Hil^{(I)})$. In the corollary below, we take $\Bw(\K_I(\fB))$ to mean the WOT sequential closure of $\K_I(\fB)$ in $B(\Hil^{(I)})$. (This is a special case of the obvious $\bw$-analogue of the spatial $C^*$-algebra tensor product, which is used below in the construction of the $\bw$-module exterior tensor product.)

Let $C_I(\fB)$ denote the collection of operators in $B(\Hil, \Hil^{(I)})$ of the form $\zeta \mapsto (b_i(\zeta))_{i \in I}$ for $b_i \in \fB$, $\zeta \in \Hil$.  This coincides with the notation in \cite[1.2.26]{BLM}, and it is shown there that $C_I(\fB)$ is a norm-closed TRO in $B(\Hil, \Hil^{(I)})$. Denote by $C_I^{\sigma}(\fB)$ the WOT sequential closure of $C_I(\fB)$ in $B(\Hil, \Hil^{(I)})$.

\begin{cor} \label{fund Me examp}
The space $C_I^{\sigma}(\fB)$ is a $\Bw(\K_I(\fB))-\fB$ $\bw$-imprivitivity bimodule.
\end{cor}

\proof
This follows immediately from Theorem~\ref{T corners are me} by noting that
\[ \left[
	\begin{matrix}
		\Bw(\K_I(\fB))	& C_I^{\sigma}(\fB) \\
		\overline{C_I^{\sigma}(\fB)} & \fB
	\end{matrix}
\right] \]
is a $\bw$-algebra in $B(\Hil^{(I)} \oplus \Hil)$ (indeed, it is the WOT sequential closure of the $\K_I(\fB)-\fB$ linking algebra of $C_I(\fB)$ in $B(\Hil^{(I)} \oplus \Hil)$) with complementary $\sigma$-full corners $\Bw(\K_I(\fB))$ and $\fB$.
\endproof

The proposition below gives a generalization and simpler proof of \cite[Theorem~3.12]{Bear}, which gives existence of the $\bw$-module completion of a $C^*$-module $X$ over a $\bw$-algebra (but it is not clear how one could deduce the facts about $\Bw(X)$ in \cite[Proposition~3.14]{Bear} from this shorter route). For a right $C^*$-module $X$ over a $C^*$-algebra $B$, denote by $\mathcal L$ the linking algebra of $X$; i.e., $\mathcal{L} := \left[ \begin{matrix} \K_B(X) & X \\ \overline X & B \end{matrix} \right]$. It is well-known that this is a $C^*$-algebra with the canonical product and involution coming from the module action and inner product on $X$, and that if $B$ is represented faithfully and nondegenerately on a Hilbert space $\Hil$, then there is a canonical faithful corner-preserving representation of $\mathcal L$ on $B((X \otimes_B \Hil) \oplus \Hil)$.

\begin{prop} \label{bw closures}
Let $X$ be a right $C^*$-module over a nondegenerate $C^*$-subalgebra $B \subseteq B(\Hil)$, and denote by $\Bw(X)$ the WOT sequential closure of $X$ in $B(\Hil, X \otimes_B \Hil)$. Then $\Bw(X)$ is a $\bw$-module over $\Bw(B)$ with inner product and module action extending those of $X$. Additionally, if $S$ is a $\sigma_{\Bw(B)}$-closed subset of $\Bw(X)$ with $X \subseteq S$, then $S = \Bw(X)$.
\end{prop}

\proof
As mentioned above, $\mathcal L$ may be viewed as a $C^*$-subalgebra of $B((X \otimes_B \Hil) \oplus \Hil)$. Taking the WOT sequential closure there, we obtain a concrete $\bw$-algebra $\Bw(\mathcal L)$, and it is elementary to argue that $\Bw(\mathcal L)$ may be identified with the space of $2 \times 2$ matrices \[ \left[ \begin{matrix} \Bw(\K_B(X)) & \Bw(X) \\ \Bw(\overline X) & \Bw(B) \end{matrix} \right], \] where all these WOT sequential closures are taken with respect to the appropriate Hilbert space representations (e.g., $\Bw(\overline X)$ is the WOT sequential closure of $\overline X$ in $B(X \otimes_B \Hil, \Hil)$). By \cite[Theorem~3.10~(2)]{Bear}, $\Bw(X)$ is a $\bw$-module over $\Bw(B)$ with inner product and module action derived from the multiplication in $\Bw(\mathcal L)$.

An easy argument using the algebra structure of $\Bw(\mathcal L)$ shows that a subset $S \subseteq \Bw(X)$ is $\sigma_{\Bw(B)}$-closed in $\Bw(X)$ if and only if $S$ is WOT sequentially closed as a subset of $B(\Hil, X \otimes_B \Hil)$, which proves the ``additionally" claim.
\endproof

Now let $X$ be a right $C^*$-module over a $C^*$-algebra $B$. As shown in \cite[Proposition~8.5.17]{BLM}, $X^{**}$ admits a canonical left $B^{**}$-action and $B^{**}$-valued inner product under which it becomes a right $W^*$-module over the $W^*$-algebra $B^{**}$, and if additionally $X$ is an $A-B$ $C^*$-imprivitivity bimodule, then $X^{**}$ is an $A^{**}-B^{**}$ $W^*$-imprivitivity bimodule.

Extending the definition of the Davies-Baire envelope of a $C^*$-algebra (see above Proposition~\ref{weak me for envs}), denote by $\Sigma(X)$ the weak* sequential closure of $X$ in $X^{**}$. As another application of Theorem~\ref{T
corners are me}, we prove in the proposition below a $\bw$-analogue of the result in the last paragraph, which essentially just notes that the structure $X^{**}$ has as a $W^*$-module over $B^{**}$ ``restricts" to $\Sigma(X)$ to make it a $\bw$-module over $\Sigma(B)$.

%

\begin{prop}
If $X$ is a right $C^*$-module over a $C^*$-algebra $B$, then $\Sigma(X)$ is canonically a right $\bw$-module over $\Sigma(B)$ whose action is the restriction of the canonical one of $B^{**}$ on $X^{**}$. If $X$ is an $A-B$ $C^*$-imprivitivity bimodule between $A$ and $B$, then $\Sigma(X)$ is a $\Sigma(A)-\Sigma(B)$ $\bw$-imprivitivity bimodule.
\end{prop}

\proof
For the first statement, let $\mathcal L := \left[ \begin{matrix} \K_B(X) & X \\ \overline{X} & B \end{matrix} \right]$. By basic functional analysis, $\mathcal L^{**} =  \left[ \begin{matrix} \K_B(X)^{**} & X^{**} \\ \overline{X}^{**} & B^{**} \end{matrix} \right],$ and by the proof of \cite[8.5.17]{BLM}, the action of $B^{**}$ on $X^{**}$ induced by $\mathcal L^{**}$ coincides with the canonical action. Routine arguments show that $\Sigma(\mathcal L) = \left[ \begin{matrix} \Sigma(\K_B(X)) & \Sigma(X) \\ \Sigma(\overline{X}) & \Sigma(B) \end{matrix} \right]$, so the result follows from \cite[Theorem~3.10~(2)]{Bear}.

The final statement in the proposition follows quickly from a similar argument and Theorem~\ref{T corners are me}. Indeed, if $\mathcal L := \left[ \begin{matrix} A & X \\ \overline{X} & B \end{matrix} \right]$, then $\Sigma(\mathcal L) = \left[ \begin{matrix} \Sigma(A) & \Sigma(X) \\ \Sigma(\overline{X}) & \Sigma(B) \end{matrix} \right]$, and $\sigma$-fullness of the 1-1 and 2-2 corners holds since $\Bw({}_{\Sigma(A)}\langle \Sigma(X)|\Sigma(X) \rangle) \supseteq \overline{{}_A \langle X|X \rangle} = A$ and $\Bw(\langle \Sigma(X)|\Sigma(X) \rangle)_{\Sigma(B)} \supseteq \overline{\langle X|X \rangle_B} = B$.
\endproof

\section{The $\Sigma^*$-module exterior tensor product}

For a right $C^*$-module $X$ over a $C^*$-algebra $B \subseteq B(\Hil)$, let $\Bw(X)$ be the WOT sequential closure of $X$ in $B(\Hil, X \otimes_B \Hil)$, so that by Proposition~\ref{bw closures}, $\Bw(X)$ is a $\bw$-algebra over $\Bw(B)$. Hence we may also view $\Bw(X)$ in $B(\Hil, \Bw(X) \otimes_{\Bw(B)} \Hil)$, but as the next lemma shows, this is really no different from viewing $\Bw(X)$ in $B(\Hil, X \otimes_B \Hil)$.

\begin{lemma} \label{lem for the lem for the lem}
If $X$ is a right $C^*$-module over a $C^*$-algebra $B \subseteq B(\Hil)$, then there is a canonical unitary $U: X \otimes_B \Hil \to \Bw(X) \otimes_{\Bw(B)} \Hil$. The spatial isomorphism between $B(\Hil, X \otimes_B \Hil)$ and $B(\Hil, \Bw(X) \otimes_{\Bw(B)} \Hil)$ induced by $U$ restricts to the identity between the canonical copies of $\Bw(X)$.
\end{lemma}

\proof
Standard arguments show that there is an isometry $X \otimes_B \Hil \to \Bw(X) \otimes_{\Bw(B)} \Hil$ determined by the rule $x \otimes \zeta \mapsto x \otimes \zeta$ on simple tensors. A calculation using the WOT sequential closure of the linking algebra of $X$ shows that \[ \langle \xi_1(h_1), \xi_2(h_2) \rangle_{X \otimes_B \Hil} = \langle h_1, \langle \xi_1| \xi_2 \rangle_{\Bw(B)} h_2 \rangle_\Hil \] for $\xi_1, \xi_2 \in \Bw(X)$ and $h_1, h_2 \in \Hil$. It follows that the rule $\xi \otimes \eta \mapsto \xi(\eta)$, for $\xi \in \Bw(X) \subseteq B(\Hil, X \otimes_B \Hil)$ and $\eta \in \Hil$, determines an isometry $\Bw(X) \otimes_{\Bw(B)} \Hil \to X \otimes_B \Hil$, which is easily seen to be the inverse of $U$.

The last statement is easy to check.
\endproof

\begin{cor} \label{totald}
The set $\{x \otimes \zeta: x \in X, \, \zeta \in \Hil \}$ is total in $\Bw(X) \otimes_{\Bw(B)} \Hil$.
\end{cor}

We now work out some of the details of the $\bw$-analogue of the \emph{exterior tensor product} of $C^*$-modules, which we will then apply to prove a fact needed in the proof of our $\bw$-analogue of the Brown-Green-Rieffel stable isomorphism theorem (mimicking the route taken to prove the stable isomorphism theorem in \cite{Lan} and \cite{RW}). To this end, let $\fA \subseteq B(\Kil)$ and $\fB \subseteq B(\Hil)$ be concrete $\bw$-algebras (although the definition works fine if $\fA$ and $\fB$ are merely $C^*$-algebras); let $\fA \otimes \fB \subseteq B(\Kil \otimes^2 \Hil)$ denote the spatial $C^*$-algebra tensor product of $\fA$ and $\fB$; and let $\fA \otimes^{\bw} \fB$ denote the $\sigma$-closure of $\fA \otimes \fB$ in $B(\Kil \otimes^2 \Hil)$, that is,
\[ \fA \otimes^{\bw} \fB := \Bw_{\Kil \otimes^2 \Hil} (\fA \otimes \fB). \]
Call $\fA \otimes^{\bw} \fB$ the \emph{$\bw$-spatial tensor product of $\fA$ and $\fB$}. (This is different from Dang's $\bw$-tensor product from \cite[B.III]{Dang}, which is the $\bw$-analogue of the maximal $C^*$-algebra tensor product.)

For the remainder of this section, take $\fX$ and $\fY$ to be right $\bw$-modules over concrete $\bw$-algebras $\fA \subseteq B(\Kil)$ and $\fB \subseteq B(\Hil)$ respectively.
Let $\fX \otimes \fY$ denote the $C^*$-module exterior tensor product of $\fX$ and $\fY$ (see e.g.\ \cite[pages 34--38]{Lan}), which is a right $C^*$-module over $\fA \otimes \fB \subseteq B(\Kil \otimes^2 \Hil)$. Denote by $\fX \otimes^{\bw} \fY$ the $\bw$-module completion of $\fX \otimes \fY$ from \cite[3.11--3.13]{Bear}, i.e.
\[ \fX \otimes^{\bw} \fY := \Bw(\fX \otimes \fY), \]
where the $\sigma$-closure is taken in $B(\Kil \otimes^2 \Hil, (\fX \otimes \fY) \otimes_{\fA \otimes \fB} (\Kil \otimes^2 \Hil))$. This is called the \emph{$\bw$-module exterior tensor product of $\fX$ and $\fY$}.


\begin{lemma} \label{lem for the lem}
\begin{enumerate}

\item[$\mathrm{(1)}$]
If $(\xi_n), \xi \in \fX \otimes^{\bw} \fY$ are such that
$\xi_n \rmto{\fA \otimes^{\bw} \fB} \xi$, then
$|\xi_n \rangle \langle \gamma | \ssto |\xi \rangle \langle \gamma |$ in $\Bw(\K_{\fA \otimes^{\bw} \fB}(\fX \otimes^{\bw} \fY))$ for all $\gamma \in \fX \otimes^{\bw} \fY$.

\item[$\mathrm{(2)}$]
If $S \in \Bw(\K_\fB(\fY))$ and $T_n \ssto T$ in $\Bw(\K_\fA(\fX))$, then $T_n \otimes S \ssto T \otimes S$ in $\Bw(\K_\fA(\fX)) \otimes^{\bw} \Bw(\K_\fB(\fY))$. 

\end{enumerate}
\end{lemma}

\proof
(1) follows by a direct application of \cite[Lemma~3.7]{Bear}.

For (2), recall that $\Bw(\K_\fA(\fX))$ and $\Bw(\K_\fB(\fY))$ can be viewed as concrete $\bw$-algebras in $B(\fX \otimes_\fA \Hil)$ and $B(\fY \otimes_\fB \Kil)$, so that $\Bw(\K_\fA(\fX)) \otimes^{\bw} \Bw(\K_\fB(\fY))$ is a concrete $\bw$-algebra in $B((\fX \otimes_\fA \Kil) \otimes^2 (\fY \otimes_\fB \Hil))$. It thus suffices to check convergence on tensors of the form $(x \otimes \zeta) \otimes (y \otimes \eta)$ since $(T_n \otimes S)$ is bounded and tensors of this form are total in $(\fX \otimes_\fA \Kil) \otimes^2 (\fY \otimes_\fB \Hil)$ by elementary arguments. The required calculation is straightforward.
\endproof

\begin{lemma} \label{ext tp dense}
The subspaces
\[ \mathrm{span}\{ |x_1 \otimes y_1 \rangle \langle x_2 \otimes y_2 | : x_1,x_2 \in \fX; \ y_1, y_2 \in \fY \} \]
and
\[ \mathrm{span}\{ |x_1 \rangle \langle x_2| \otimes |y_1 \rangle \langle y_2| : x_1,x_2 \in \fX; \ y_1, y_2 \in \fY \} \]
are $\sigma$-dense in $\Bw(\K_{\fA \otimes^{\bw} \fB}(\fX \otimes^{\bw} \fY))$ and $\Bw(\K_\fA(\fX)) \otimes^{\bw} \Bw(\K_\fB(\fY))$ respectively.
\end{lemma}

\proof
Call the displayed sets $V_1$ and $V_2$ respectively. These are $*$-algebras, so that $\Bw(V_1)$ and $\Bw(V_2)$ are $\bw$-algebras. We have the evident inclusions
\[
\Bw(V_1) \subseteq \Bw(\K_{\fA \otimes^{\bw} \fB}(\fX \otimes^{\bw} \fY)) \]
\[
\Bw(V_2) \subseteq \Bw(\K_\fA(\fX)) \otimes^{\bw} \Bw(\K_\fB(\fY)).
\]
So to prove the lemma, it suffices to show $\K_{\fA \otimes^{\bw} \fB}(\fX \otimes^{\bw} \fY) \subseteq \Bw(V_1)$ and $\Bw(\K_\fA(\fX)) \otimes \Bw(\K_\fB(\fY)) \subseteq \Bw(V_2)$, which follow if we can show
\begin{equation} \label{L23 1}
| \xi \rangle \langle \gamma | \in \Bw(V_1) \text{ for all } \xi, \gamma \in \fX \otimes^{\bw} \fY
\end{equation}
and
\begin{equation} \label{L23 2}
T \otimes S \in \Bw(V_2) \text{ for all } S \in \Bw(\K_\fA(\fX)) \text{ and } T \in \Bw(\K_\fB(\fY)).
\end{equation}
For (1), fix $\gamma \in \fX \otimes \fY$, and let $\mathscr S = \{ \xi \in \fX \otimes^{\bw} \fY : |\xi \rangle \langle \gamma| \in \Bw(V_1)\}$. Clearly $\fX \otimes \fY \subseteq \mathscr S$, and it follows by Lemma~\ref{lem for the lem} (1) that $\mathscr S$ is $\sigma_{\fA \otimes^{\bw} \fB}$-closed in $\fX \otimes^{\bw} \fY$. By the ``additionally" statement in Proposition~\ref{bw closures}, we get that $\mathscr S =  \fX \otimes^{\bw} \fY$. A similar argument using this proves (1).

A similar argument to that in the first paragraph, using the second part of Lemma~\ref{lem for the lem}, gives (2) here.
\endproof

\begin{thm} [cf.\ \cite{Lan}, end of page 37] \label{spatial tp thm}
There is a canonical $\bw$-isomorphism \[ \Bw(\K_{\fA \otimes^{\bw} \fB}(\fX \otimes^{\bw} \fY)) \cong \Bw(\K_\fA(\fX)) \otimes^{\bw} \Bw(\K_\fB(\fY)). \]
\end{thm}

\proof
Recall as above that the two $\bw$-algebras in the conclusion of the claim can be viewed as concrete $\bw$-algebras as follows:
\begin{align*}
\Bw(\K_{\fA \otimes^{\bw} \fB}(\fX \otimes^{\bw} \fY)) &\subseteq B((\fX \otimes^{\bw} \fY) \otimes_{\fA \otimes^{\bw} \fB} (\Kil \otimes^2 \Hil)), \\
\Bw(\K_\fA(\fX)) \otimes^{\bw} \Bw(\K_\fB(\fY)) & \subseteq B((\fX \otimes_\fA \Kil) \otimes^2 (\fY \otimes_\fB \Hil)).
\end{align*}
For convenience, label the Hilbert spaces $\mathcal M := (\fX \otimes^{\bw} \fY) \otimes_{\fA \otimes^{\bw} \fB} (\Kil \otimes^2 \Hil)$ and $\mathcal N := (\fX \otimes_\fA \Kil) \otimes^2 (\fY \otimes_\fB \Hil)$, and label the $\bw$-modules $\mathfrak M := \Bw(\K_{\fA \otimes^{\bw} \fB}(\fX \otimes^{\bw} \fY))$ and $\mathfrak N := \Bw(\K_\fA(\fX)) \otimes^{\bw} \Bw(\K_\fB(\fY))$. We will show that there is a canonical unitary $U: \mathcal M \to \mathcal N$ which implements a spatial isomorphism between $\mathfrak M$ and $\mathfrak N$. To see the existence of such a unitary, first define $U$ on simple tensors of simple tensors in $\mathcal M$ by
\[ U( (x \otimes y) \otimes (\zeta \otimes \eta))= (x \otimes \zeta) \otimes (y \otimes \eta) \]
for $x \in \fX$, $y \in \fY$, $\zeta \in \Kil$, and $\eta \in \Hil$. It is then straightforward using the definitions of the inner products in these various $C^*$-modules and Hilbert spaces to see that $U$ extends to an isometry on the span of elements of the form $(x \otimes y) \otimes (\zeta \otimes \eta)$ in $\mathcal M$.
Since the set of elements of the form $ (x \otimes y) \otimes (\zeta \otimes \eta)$ (resp.\ $(x \otimes \zeta) \otimes (y \otimes \eta)$)  is total in $\mathcal M$ (resp.\ $\mathcal N$) by Corollary~\ref{totald} (resp.\ by elementary principles), it follows that $U$ extends to a unitary $U: \mathcal M \to \mathcal N$.

Let $\varphi : B(\mathcal M) \to B(\mathcal N)$ denote the isomorphism $m \mapsto UmU^*$. A straightforward calculation gives the formula
\[ \varphi(|x_1 \otimes y_1 \rangle \langle x_2 \otimes y_2 |) = |x_1 \rangle \langle x_2| \otimes |y_1 \rangle \langle y_2| \]
for all $x_1,x_2 \in \fX$ and $y_1, y_2 \in \fY$, showing that $\varphi$ restricts to an isomorphism from
\[ \mathrm{span}\{ |x_1 \otimes y_1 \rangle \langle x_2 \otimes y_2 | : x_1,x_2 \in \fX; \ y_1, y_2 \in \fY \} \]
onto
\[ \mathrm{span}\{ |x_1 \rangle \langle x_2| \otimes |y_1 \rangle \langle y_2| : x_1,x_2 \in \fX; \ y_1, y_2 \in \fY \}. \]
It is easy to check that if $C$ is a $C^*$-subalgebra of $B(\Hil_1)$ and $V:\Hil_1 \to \Hil_2$ is a unitary implementing an isomorphism $\psi:B(\Hil_1) \to B(\Hil_2)$, then $\psi(\Bw(C)) = \Bw(\psi(C))$. Applying this fact to $\varphi$ and using Lemma~\ref{ext tp dense}, we get that $\varphi$ restricts to a $\bw$-isomorphism from $\Bw(\K_{\fA \otimes^{\bw} \fB}(\fX \otimes^{\bw} \fY))$ onto $\Bw(\K_\fA(\fX)) \otimes^{\bw} \Bw(\K_\fB(\fY))$.
\endproof

\section{A Brown-Green-Rieffel stable isomorphism theorem for $\Sigma^*$-modules}

We conclude by proving a $\bw$-analogue of the Brown-Green-Rieffel stable isomorphism theorem, which asserts that $\sigma$-unital $C^*$-algebras are strongly $C^*$-Morita equivalent if and only if they are stably isomorphic. (See any of the texts mentioned in Section~2 for a proof, and see \cite[8.5.31]{BLM} for the $W^*$-version.) Unfortunately, we have not been able to find a condition on the algebras in our setting as elegant as the condition of being $\sigma$-unital in the $C^*$-setting (the basic problem is that we need a condition on the algebras guaranteeing that some $\bw$-imprivitivity bimodule is $\bw$-countably generated on both sides). Our method of proof is a translation to the $\bw$-setting of the proofs given in \cite[8.2.7]{BLM} (for $C^*$-algebras) and \cite[8.5.31]{BLM} (for $W^*$-algebras). We first recall some definitions and results from \cite[Section~4]{Bear} that we will need.

\begin{defn} [\cite{Bear}, Definition~4.1]
A right (resp.\ left) $\bw$-module $\fX$ over a $\bw$-algebra $\fA \subseteq B(\Hil)$ is \emph{$\Sigma^*_\fA$-countably generated} (resp.\ ${}_\fA \Sigma^*$-countably generated) if there is a countable set $\{x_i\}_{i=1}^\infty \subseteq \fX$ such that $\{\sum_{i=1}^N x_i a_i : a_i \in \fA, N \in \N\}$ (resp.\ $\{\sum_{i=1}^N a_i x_i : a_i \in \fA, N \in \N\}$) is WOT sequentially dense in $\fX \subseteq B(\Hil, \fX \otimes_\fA \Hil)$ (resp.\ in $\fX \subseteq B(\overline \fX \otimes_\fA \Hil, \Hil)$).
\end{defn}

\begin{defn} [cf.\ \cite{Bear}, Lemma~4.14 and above]
If $\fX$ is a $\bw$-module over $\fB$, denote by $C^w(\fX)$ the $\bw$-module with underlying vector space
\[ \{(x_n) \in \prod_n \fX : \sum_n \langle x_n|x_n \rangle \text{ is } \sigma \text{-convergent in } \fB\}, \]
inner product $\langle (x_n)|(y_n) \rangle := \sum_n \langle x_n| y_n \rangle$, and entrywise $\fB$-module action. It is not hard to see from the definitions that if $\fX = \fB$, this coincides with $C_\mathbb{N}^{\sigma}(\fB)$ from above Corollary~\ref{fund Me examp} (more precisely, there is a canonical unitary between $C^w(\fB)$ and $C_\mathbb{N}^{\sigma}(\fB)$ induced by the canonical unitary between $\Hil^{(\mathbb N)}$ and $C^w(\fB) \otimes_\fB \Hil$).
\end{defn}

\begin{defn} [cf.\ \cite{Bear}, Definition~4.23]
For a left $\bw$-module $\fX$ over a $\bw$-algebra $\fC$, a countable subset $\{x_k\}$ of $\fX$ is a called a \emph{left weak quasibasis} for $\fX$ if for any $x \in \fX,$ the sequence of finite sums ($\sum_{k=1}^n \langle x| x_k \rangle x_k)_n$ ${}_\fC \sigma$-converges to $x$. A similar definition holds for \emph{right weak quasibases}.
\end{defn}

\begin{prop} [cf.\ \cite{Bear}, Proposition~4.7]
 \label{wcg}
Let $\fX$ be a right $\bw$-module over a concrete $\bw$-algebra $\fB \subseteq B(\Hil).$ Then $\fX$ is $\swb$countably generated if and only if there is a sequence $(e_n)$ in $\K_\fB(\fX)$ such that $e_n \ssto I$ in $\B_\fB(\fX)$ (recall the $\bw$-algebra structure on $\B_\fB(\fX)$ mentioned above Proposition~\ref{Me equiv lr}).
\end{prop}

\begin{thm} [\cite{Bear}, Theorem~4.19]
\label{stabilization}
If $\fB$ is a $\bw$-algebra and $\fX$ is a $\swb$countably generated right $\bw$-module over $\fB,$ then $\fX \oplus C^w(\fB) \cong C^w(\fB)$ unitarily.
\end{thm}

\begin{prop} [cf.\ \cite{Bear}, Proposition~4.21, noting that part of the assumption in (3) is in fact not necessary]
\label{orth compl conds}
Let $\fX$ be a submodule of a $\bw$-module $\fY$ over a $\bw$-algebra $\fB$. If $\fX$ meets the following requirements:
\begin{enumerate}
\item[$\mathrm{(1)}$] $\fX$ is $\sigma_\fB$-closed in $\fY$;
\item[$\mathrm{(2)}$] $\Bw(\K_\fB(\fX)) = \B_\fB(\fX)$;
\end{enumerate}
then $\fX$ is orthogonally complemented in $\fY$.
\end{prop}

\begin{thm} [cf.\ \cite{Bear}, Theorem~4.26]
\label{lwq}
Let $\fC \subseteq B(\Hil)$ be a unital $\bw$-algebra, and let $\fX$ be a left $\bw$-module over $\fC.$ Then $\fX$ is ${}_\fC \bw$-countably generated if and only if $\fX$ has a left weak quasibasis.
\end{thm}

\begin{prop} \label{colsp is etp}
If $\fX$ is a $\bw$-module over a $\bw$-algebra $\fB$, then there is a canonical unitary between $C^w(\fX)$ and the $\bw$-module exterior tensor product $\ell^2 \otimes^{\bw} \fX$.
\end{prop}

\proof
Fix a faithful $\bw$-representation $\fB \hookrightarrow B(\Hil)$. By the proof of \cite[Lemma~4.14]{Bear}, $C^w(\fX)$ is unitarily equivalent to $\mathscr C:= \{ T \in B(\Hil, (\fX \otimes_\fB \Hil)^{(\mathbb N)}) : P_n T \in \fX\}$, where $P_n : (\fX \otimes_\fB \Hil)^{(\mathbb N)} \to \fX \otimes_\fB \Hil$ is the projection onto the $n^{\mathrm{th}}$ coordinate. This unitary restricts to a unitary between the closure of \[ \{T \in \mathscr C : P_nT = 0 \text{ for all but finitely many } n\} \] and $C(\fX)$ (this is the $C^*$-module direct sum of countably many copies of $\fX$). It is easy to argue from this that $C(\fX)$ is WOT sequentially dense in $C^w(\fX)$ when these are viewed in $B(\Hil, (\fX \otimes_B \Hil)^{(\mathbb N)})$. By basic properties of the interior tensor product of $C^*$-module (see \cite[8.2.12~(3)]{BLM}), the Hilbert spaces $(\fX \otimes_\fB \Hil)^{(\mathbb N)}$ and $C(X) \otimes_\fB \Hil$ are unitarily equivalent, and it is easy to check that the induced isomorphism between $B(\Hil, (\fX \otimes_\fB \Hil)^{(\mathbb N)})$ and $B(\Hil, C(\fX) \otimes_\fB \Hil)$ restricts to the identity on the copies of $C(\fX)$ in each of these. Additionally, it is a well-known $C^*$-module fact (see \cite[pg.~35]{Lan}) that $C(\fX)$ is unitarily equivalent to the exterior tensor product $\ell^2 \otimes \fX$. Putting all these together yields the following chain of canonical unitaries:
\begin{align*}
C^w(\fX) &\cong \Bw_{B(\Hil, (\fX \otimes_\fB \Hil)^{(\mathbb N)})}(C(X)) \\
&\cong \Bw_{B(\Hil,C(X) \otimes_\fB \Hil)}(C(X))\\
&\cong \Bw_{B(\Hil,(\ell^2 \otimes \fX) \otimes_\fB \Hil)}(\ell^2 \otimes \fX)\\
&= \ell^2 \otimes^{\bw} \fX.
\end{align*}
\endproof

\begin{lemma} \label{swindle}
If $\fX$ and $\fY$ are $\bw$-modules over $\fB$, then $C^w(C^w(\fX)) \cong C^w(\fX)$ and $C^w(\fX \oplus \fY) \cong C^w(\fX) \oplus C^w(\fX \oplus \fY)$ unitarily.
\end{lemma}

\proof
By the proof of \cite[Lemma~4.14]{Bear}, if $\fB \subseteq B(\Hil)$ is a faithful $\bw$-representation, then
\[ C^w(\fX) \cong \{T \in B(\Hil,(\fX \otimes_\fB \Hil)^{(\mathbb N)}) : P_n T \in \fX \text{ for all } n \in \mathbb N\} \]
and
\[ C^w(C^w(\fX)) \cong \{T \in B(\Hil,(C^w(\fX) \otimes_\fB \Hil)^{(\mathbb N)}) : P_n' T \in C^w(\fX) \text{ for all } n \in \mathbb N\} \]
unitarily (here, $P_n: \Hil \to \fX \otimes_\fB \Hil$ and $P_n': \Hil \to C^w(\fX) \otimes_\fB \Hil$ denote the canonical projections onto the $n^{\text{th}}$ summands). Identifying $C^w(\fX) \otimes_\fB \Hil$ with $(\fX \otimes_\fB \Hil)^{(\mathbb N)}$ via the canonical unitary, it follows that $(C^w(\fX) \otimes_\fB \Hil)^{(\mathbb N)} \cong ((\fX \otimes_\fB \Hil)^{(\mathbb N)})^{(\mathbb N)}$. By properties of Hilbert space direct sums, $((\fX \otimes_\fB \Hil)^{(\mathbb N)})^{(\mathbb N)} \cong (\fX \otimes_\fB \Hil)^{(\mathbb N \times \mathbb N)} \cong (\fX \otimes_\fB \Hil)^{(\mathbb N)}$. The latter unitary is not canonical, but for any choice of such a unitary, it readily follows that the unitary $(C^w(\fX) \otimes_\fB \Hil)^{(\mathbb N)} \cong (\fX \otimes_\fB \Hil)^{(\mathbb N)}$ induces a $C^*$-module unitary between the sets displayed above.

For the second claim, one can easily check that $C^w(\fX \oplus \fY) \cong C^w(\fX) \oplus C^w(\fY)$ via the unitary $(x_n,y_n) \mapsto ((x_n),(y_n))$, and that $C^w(\fX) \cong C^w(\fX) \oplus C^w(\fX)$ via the unitary $(x_n) \mapsto (x_{2n-1},x_{2n})$. So
\begin{align*}
C^w(\fX \oplus \fY) &\cong C^w(\fX) \oplus C^w(\fY) \\
 &\cong C^w(\fX) \oplus C^w(\fX) \oplus C^w(\fY) \\
&\cong C^w(\fX) \oplus C^w(\fX \oplus \fY).
\end{align*}
\endproof

The next result is analogous to the $C^*$-result in \cite[8.2.6~(3)]{BLM} and one of the facts shown in the proof of the $W^*$-module version of Kasparov's stabilization theorem presented in \cite[8.5.28]{BLM}.

\begin{prop} \label{other side}
If $\fX$ is a $\sigma$-full right $\bw$-module over a $\bw$-algebra $\fB \subseteq B(\Hil)$ such that $\Bw(\K_\fB(\fX))$ is unital (i.e., $\Bw(\K_\fB(\fX)) = \B_\fB(\fX)$) and $\fX$ is ${}_{\Bw(\K_\fB(\fX))} \bw$-countably generated, then $C^w(\fX) \oplus C^w(\fB) \cong C^w(\fX)$ unitarily.
\end{prop}

\proof (cf.\ \cite[proof of Corollary~8.2.6]{BLM}.)
The hypotheses allow us to invoke Theorem~\ref{lwq} on the left $\bw$-module $\fX$ over $\Bw(\K_\fB(\fX))$ to obtain a left weak quasibasis $\{x_k\}$. For any $x \in \fX$, we have
\[ x \sum_{k=1}^n \langle x_k| x_k \rangle = \sum_{k=1}^n | x \rangle \langle x_k | x_k \lmto{\Bw(\K_\fB(\fX))} x,\]
so that $x \sum_{k=1}^n \langle x_k|x_k \rangle \rmto{\fB} x$ by \cite[Lemma~3.7]{Bear}. By the ``right version" of Lemma~\ref{L WOTconv} in the current paper, $\sum_{k=1}^n \langle x_k|x_k \rangle \ssto 1_{M(\fB)}$. In particular, $\fB$ is unital and $(x_k) \in C^w(\fX)$.


Now define a map $\varphi: \fB \to C^w(\fX)$ by $\varphi(b) = (x_k b)$ for $b \in \fB$. By the calculation $\sum_{k=1}^n \langle x_k b| x_k b \rangle = b^* \sum_{k=1}^n \langle x_k|x_k \rangle b \leq b^*b$, we have that $\varphi$ does indeed map into $C^w(\fB)$, and it is easy to see that $\varphi$ is a $\fB$-module map, so that $\text{Ran}(\varphi)$ is a $\fB$-submodule of $C^w(\fB)$. We now check the conditions in Proposition~\ref{orth compl conds} for the submodule $\text{Ran}(\varphi) \subseteq C^w(\fB)$. For condition (1), we need to show that if we have a sequence $(b_n)$ from $\fB$ such that $\langle \varphi(b_n) | (y_k) \rangle$ is $\sigma$-convergent for all $(y_k) \in C^w(\fB)$, then there is a $b \in \fB$ such that $\langle \varphi(b_n)| (y_k) \rangle \ssto \langle \varphi(b)| (y_k) \rangle$ for all $(y_k) \in C^w(\fB)$. Suppose then that we have such a sequence $(b_n)$ as in the former. Then the sequence with terms
\[ b_n^* = b_n^* \sum_k \langle x_k|x_k \rangle = \langle (x_k b_n)|(x_k) \rangle = \langle \varphi(b_n)| \varphi(1) \rangle \]
$\sigma$-converges in $\fB$, so that $b_n \ssto b$ for some $b \in \fB$. Then for any $(y_k) \in C^w(\fX)$,
\[ \langle \varphi(b_n)| (y_k) \rangle = \langle (x_k b_n)|(y_k) \rangle = b_n^* \sum_k \langle x_k|y_k \rangle \ssto b^* \sum_k \langle x_k|y_k \rangle = \langle \varphi(b)| (y_k) \rangle. \]
To see that $\varphi(\fB)$ meets condition (2) from Proposition~\ref{orth compl conds}; i.e., that the identity map $I_{\varphi(\fB)} \in \B_\fB(\varphi(\fB))$ is in $\Bw(\K_\fB(\varphi(\fB)))$, let $b \in \fB$ and compute
\[ |\varphi(1) \rangle \langle \varphi(1) |(\varphi(b)) = \varphi(1) \langle (x_k) |(x_k b) \rangle = \varphi(1) \sum_k \langle x_k|x_k \rangle b =\varphi(1)b = \varphi(b). \]
Hence $I_{\varphi(\fB)} = |\varphi(1) \rangle \langle \varphi(1) | \in \K_\fB(\varphi(\fB))$. So we may apply Proposition~\ref{orth compl conds} to conclude that $\varphi(\fB)$ is orthogonally complemented in $C^w(\fB)$.

The calculation
\[ \langle \varphi(b) | \varphi(c) \rangle = \langle (bx_k)|(cx_k) \rangle = b^* \left(\sum_k \langle x_k|x_k \rangle \right) c = b^*c \]
shows that $\varphi: \fB \to \varphi(\fB)$ is a unitary map, so $C^w(\fX) \cong \fB \oplus \mathfrak W$ for some $\bw$-module $\mathfrak W$ ($\mathfrak W$ is a $\bw$-module by \cite[Proposition~4.21]{Bear}). By Lemma~\ref{swindle}, we have
\begin{equation*}
\begin{split}
C^w(\fX) \cong C^w(C^w(\fX)) \cong C^w(\fB \oplus \mathfrak W) \cong C^w(\fB) \oplus C^w(\fB \oplus \mathfrak W) \\ \cong C^w(\fB) \oplus C^w(C^w(\fX)) \cong C^w(\fB) \oplus C^w(\fX).
\end{split}
\end{equation*}
\endproof


\begin{defn}
Say that a unital $\bw$-algebra $\fB$ has \emph{Property~(D)} if for every closed ideal $I$ in $\fB$ such that $\Bw(I)=\fB=M(I)$, there is a sequence $(e_n)$ in $I$ with $e_n \ssto 1$.
\end{defn}

\begin{note} \label{prop D note}
Clearly every simple unital $\bw$-algebra (in particular, every von Neumann algebra factor that is type I$_n$, type II$_1$, or countably decomposable and type III) has Property~(D). Every infinite type I von Neumann algebra factor also has Property~(D). Indeed, if $\Hil$ is infinite dimensional and nonseparable, then $B(\Hil)$ has no WOT sequentially dense ideals; if $\Hil$ is infinite dimensional and separable, then the unique closed ideal $\mathbb K \subseteq B(\Hil)$ meets the conditions and conclusion in the definition.

\end{note}

For a $\bw$-algebra $\fB \subseteq B(\Hil)$, note that the WOT sequential closure of $\mathbb K \otimes \fB$ in $B(\ell^2 \otimes^2 \Hil)$, is all of $\mathbb{M}(\fB)$ (where the latter is the space of infinite matrices over $\fB$ indexed by $\mathbb N$ with uniformly bounded finite submatrices---see \cite[1.2.26]{BLM}). Indeed, take $A \in \mathbb{M}(\fB)$ viewed as a matrix with entries in $\fB$, and let $A_n \in \mathbb K \otimes \fB$ be the finite submatrix of $A$ supported in the entries $1,\dots,n$. An easy argument (using finitely supported columns in $\Hil^{(\mathbb N)}$ and Lemma~\ref{simplem} if you like) shows that $A_n \xrightarrow{WOT} A$. Since $\mathbb K \otimes \fB \subseteq \Bw(\mathbb K) \otimes \fB \subseteq B(\ell^2 \otimes^2 \Hil)$, it follows that the WOT sequential closure of $\Bw(\mathbb K) \otimes \fB$; i.e., $\Bw(\mathbb K) \otimes^{\bw} \fB$ in the notation introduced above Lemma~\ref{lem for the lem}, coincides with $\mathbb{M}(\fB)$.

\begin{thm}
Let $\fA$ and $\fB$ be two unital $\bw$-algebras with Property~(D). Then $\fA$ and $\fB$ are $\bw$-Morita equivalent if and only if $\mathbb{M}(\fA) \cong \mathbb{M}(\fB)$ $\bw$-isomorphically.
\end{thm}

\proof
$(\implies)$ Let $\fX$ be an $\fA-\fB$ $\bw$-imprivitivity bimodule. Then $\fA \cong \Bw(\K_\fB(\fX))$ $\bw$-isomorphically by Theorem~\ref{Me equiv}, so by Property~(D), there is a sequence $(e_n)$ in $\K_\fB(\fX)$ with $e_n \ssto I$ in $\B_\fB(\fX) = \Bw(\K_\fB(\fX))$. By Proposition~\ref{wcg}, $\fX$ is $\bw_\fB$-countably generated. It is similarly shown that $\fX$ is $\bw_\fA$-countably generated. By Theorem~\ref{stabilization}, $C^w(\fB) \cong \fX \oplus C^w(\fB)$, so that by Lemma~\ref{swindle},
\[ C^w(\fB) \cong C^w(\fX \oplus C^w(\fB)) \cong C^w(\fX) \oplus C^w(\fB). \]
By Proposition~\ref{other side}, we have $C^w(\fX) \oplus C^w(\fB) \cong C^w(\fX)$. Hence $C^w(\fX) \cong C^w(\fB)$ unitarily. By symmetry and the other-handed versions of everything above, also $C^w(\fX) \cong C^w(\fA)$ unitarily. Hence
\begin{align*}
\mathbb M(\fB)
\cong \Bw(\mathbb K) \otimes^{\bw} \fB
\cong \Bw(\K_\fB(\ell^2 \otimes^{\bw} \fB))
\cong \Bw(\K_\fB(C^w(\fB))) 
\cong \Bw(\K_\fB(C^w(\fX))) \\
\cong \Bw(\K_\fB(\ell^2 \otimes^{\bw} \fX))
\cong \Bw(\mathbb K) \otimes^{\bw} \Bw(\K_\fB(\fX))
\cong \Bw(\K) \otimes^{\bw} \fA
\cong \mathbb M(\fA),
\end{align*}
where we have used the observation above the statement of the current theorem for the first and last isomorphisms, Theorem~\ref{spatial tp thm} for the second and sixth isomorphisms (using the $\bw$-module $\ell^2$ over $\Bw(\mathbb K)$ and the fact that $\Bw(\K_\fB(\fB)) \cong \fB$ $\bw$-isomorphically) and Proposition~\ref{colsp is etp} for the third and fifth isomorphisms. Since these are all $\bw$-isomorphisms, the claim is proved.

$(\impliedby)$ This direction follows immediately from our results that $\bw$-Morita equivalence is an equivalence relation coarser than $\bw$-isomorphism (Theorem~\ref{T er}), that $\fA$ is $\bw$-Morita equivalent to $\Bw(\K_I(\fA))$ for all $I$ (Corollary~\ref{fund Me examp}), and the fact mentioned above the theorem that $\Bw(\K(\fA)) = \mathbb M(\fA)$ (and similarly for $\fB$).
\endproof

%
%





\begin{thebibliography}{99}

\bibitem{Bear}
C. A. Bearden, \textit{Hilbert $C^*$-modules over $\bw$-algebras}, Studia Math. {\bf 235} (2016), 269--304.

\bibitem{Beer}
W. Beer, \textit{On Morita equivalence of nuclear $C^*$-algebras}, J. Pure Appl. Algebra \textbf{26} (1982), 249--267.

\bibitem{Bl}
B. Blackadar, \textit{Operator algebras---theory of $C^*$-algebras and von Neumann algebras}, Encyclopaedia Math. Sci. 122, Springer-Verlag, Berlin 2006.



\bibitem{Ble97b}
D. P. Blecher, \textit{On selfdual Hilbert modules, Operator algebras and their applications}, Fields Inst. Commun. \textbf{13} (1997), 65--80.
%
%
%
\bibitem{BLM}
D. P. Blecher and C. Le Merdy, \textit{Operator algebras and their modules---an operator space approach}, Oxford Univ. Press, Oxford, 2004.

\bibitem{Dang}
N. N. Dang, \textit{$\Sigma^*$-algebras, probabilit\'{e}s non commutatives et applications}, M\'{e}moires de la S. M. F. \textbf{35} (1973), 145--189.

\bibitem{Dav68}
E. B. Davies, \textit{On the Borel structure of $C^*$-algebras}, Comm. Math. Phys. \textbf{8} (1968), 147--163.
%
\bibitem{Dav69}
E. B. Davies, \textit{The structure of $\Sigma^*$-algebras}, Quart. J. Math. Oxford. (2) \textbf{20} (1969), 351--366.
%

\bibitem{ERRS}
S. Eilers, G. Restorff, E. Ruiz, and A. P. W. S\o rensen, \textit{The complete classification of unital graph $C^*$-algebras: Geometric and strong}, preprint (arXiv: 1611.07120).
%
%
%
\bibitem{Lan}
E. C. Lance, \textit{Hilbert $C^*$-modules: A toolkit for operator algebraists}, London Math. Soc. Lecture Note Series, vol. 210, Cambridge Univ. Press, Cambridge, 1994.
%
%

\bibitem{SML}
S. MacLane, \textit{Categories for the Working Mathematician}, Springer-Verlag, New York, 1998.


\bibitem{Munkres}
J. R. Munkres, \textit{Topology}, $2^{\mathrm{nd}}$ edition, Prentice-Hall, Upper Saddle River (2000).
%
%
%
\bibitem{Ped NOW}
G. K. Pedersen, \textit{Analysis now}, Graduate Texts in Math., vol. 118, Springer-Verlag, New York, 1989. (Corrected second printing, 1995.)
%

\bibitem{RW}
I. Raeburn and D. P. Williams, \textit{Morita equivalence and continuous-trace $C^*$-algebras}, Mathematical Surveys and Monographs, vol. 60, American Mathematical Society, Providence, RI, 1998.

\bibitem{Ror}
M. R\o rdam, \textit{Classification of Nuclear, Simple $C^*$-algebras}, Encyclopaedia of Mathematical Sciences, vol. 126, Springer, Berlin, 2002.

\bibitem{Rie}
M. Rieffel, \textit{Morita equivalence for $C^*$-algebras and $W^*$-algebras}, J. Pure Appl. Algebra \textbf{5} (1974), 51--96.

\bibitem{Rie2}
M. Rieffel, \textit{Induced representations of $C^*$-algebras}, Adv. Math. \textbf{13} (1974), 176--257.

%
%

\end{thebibliography}
\end{document}